\documentclass[10pt,twoside]{article}
\usepackage[latin1]{inputenc}
\usepackage{amsmath}
\usepackage{graphicx}
\usepackage{yhmath}
\usepackage[table]{xcolor}
\usepackage{mathrsfs} 
\usepackage{amssymb}
\usepackage{url}
\usepackage{makecell}
\usepackage{arydshln}
\usepackage{multirow}
\usepackage{arydshln}
\usepackage{mathtools}
\usepackage{stmaryrd}  

\input epsf
\def\limiten{\renewcommand{\arraystretch}{0.5}
\begin{array}[t]{c}\stackrel{}{\longrightarrow} \\
{\scriptstyle n\rightarrow
\infty}\end{array}\renewcommand{\arraystretch}{1}}

\def\limitepsn{\renewcommand{\arraystretch}{0.5}
\begin{array}[t]{c}\stackrel{a.s.}{\longrightarrow} \\
{\scriptstyle n \rightarrow
\infty}\end{array}\renewcommand{\arraystretch}{1}}

\def\limiteloin{\renewcommand{\arraystretch}{0.5}
\begin{array}[t]{c}\stackrel{{\cal D}}{\longrightarrow} \\
{\scriptstyle n\rightarrow
\infty}\end{array}\renewcommand{\arraystretch}{1}}

\def\limiteproban{\renewcommand{\arraystretch}{0.5}
\begin{array}[t]{c}\stackrel{{\cal P}}{\longrightarrow} \\
{\scriptstyle n\rightarrow
\infty}\end{array}\renewcommand{\arraystretch}{1}}

\def\equalpsn{\renewcommand{\arraystretch}{0.5}
\begin{array}[t]{c}\stackrel{a.s.}{=} \\
\end{array}\renewcommand{\arraystretch}{1}}

\numberwithin{equation}{section}

\newtheorem{thm}{Theorem}[section]

\newtheorem{lem}[thm]{Lemma}

\newtheorem{prop}[thm]{Proposition}

\newcommand{\E}{\ensuremath{\mathbb{E}}}
\newcommand{\R}{\ensuremath{\mathbb{R}}}
\newcommand{\Z}{\ensuremath{\mathbb{Z}}}

\newcommand{\N}{\ensuremath{\mathbb{N}}}

\newcommand{\prob}{\ensuremath{\mathbb{P}}}

\definecolor{grisclair}{gray}{0.9}

\renewcommand{\arraystretch}{.8}

\begin{document}
\title{\bf A general procedure for change-point detection in multivariate time series}
 \maketitle \vspace{-1.0cm}
  \begin{center}
   Mamadou Lamine DIOP \footnote{Supported by
   the MME-DII center of excellence (ANR-11-LABEX-0023-01) 
   } 
   and 
    William KENGNE \footnote{Developed within the ANR BREAKRISK: ANR-17-CE26-0001-01 and the  CY Initiative of Excellence (grant "Investissements d'Avenir" ANR-16-IDEX-0008), Project "EcoDep" PSI-AAP2020-0000000013} 
 \end{center}

  \begin{center}
  { \it 
 THEMA, CY Cergy Paris Université, 33 Boulevard du Port, 95011 Cergy-Pontoise Cedex, France.\\
  E-mail: mamadou-lamine.diop@u-cergy.fr ; william.kengne@u-cergy.fr  \\
  }
\end{center}

 \pagestyle{myheadings}
 \markboth{ A general procedure for change-point detection in multivariate time series}{Kengne and Diop}

~~\\
\textbf{Abstract}: We consider the change-point detection in multivariate continuous and integer valued time series. We propose a Wald-type statistic based on the estimator performed by a general contrast function; which can be constructed from the likelihood, a quasi-likelihood, a least squares method, etc. 
 Sufficient conditions are provided to ensure that the statistic convergences to a well-known distribution under the null hypothesis (of no change) and  diverges to infinity under the alternative; which establishes the consistency of the procedure. Some examples are detailed to illustrate the scope of application of the proposed procedure. 
 Simulation experiments are conducted to illustrate the asymptotic results.\\


 \section{Introduction}\label{sec_intro}
 %
%
 
%
Since \nocite{Page1955} Page (1955), the change-point problem  has been widely studied. Several approaches and procedures have been developed for continuous and integer valued, univariate and multivariate processes.  \\
 Consider observations $(Y_{1},\cdots,Y_{n})$, generated from a multivariate continuous or integer valued process $Y=\{Y_{t},\,t\in \Z \}$. These observations depend on a parameter $\theta^* \in \Theta \subset \R^d $ ($d \in \N$) which may change over time.  
More precisely, consider the following test hypotheses:

\begin{enumerate}
    \item [ H$_0$:] $(Y_1,\cdots,Y_n)$  is a trajectory of the process $Y=\{Y_{t},\,t\in \Z \}$ which depends on $\theta^*$.
    
     \item[ H$_1$:] There exists $((\theta^{*}_1,\theta^{*}_2),t^{*}) \in \Theta^{2}\times \{2,3,\cdots, n-1 \}$ (with $\theta^{*}_1 \neq \theta^{*}_2$) such that,
     $(Y_1,\cdots,Y_{t^{*}})$ is a trajectory of a process $Y^{(1)} = \{Y^{(1)}_{t},\, t \in \Z\}$ that depends on $\theta^*_1$ and  $(Y_{t^{*}+1},\cdots,Y_n)$ is a trajectory of a process $Y^{(2)} = \{Y^{(2)}_{t}, \,t \in \Z\}$ that depends on $\theta^*_2$. 
\end{enumerate}
Note that, under H$_1$, $(Y_1,\cdots,Y_n)$ is a trajectory of the process $\{(Y^{(1)}_{t})_{t \leq t^*}, (Y^{(2)}_{t})_{t > t^*} \}$ which depends on $\theta^*_1$ and $\theta^*_2$.
In the whole paper, it is assumed that $\Theta$ is a fixed compact subset of $\R^d$ ($d \in \N$). 

\medskip

This test for change-point detection is often addressed with a Wald-type statistic based on the likelihood, quasi likelihood, conditional least-squares or density power divergence estimator.
Likelihood estimate-based procedure has been proposed for  continuous and integer valued time series; see for instance, \nocite{Lee2004cusum} Lee and Lee (2004), \nocite{Kang2014} Kang and Lee (2014), \nocite{Doukhan2015} Doukhan and Kengne (2015), \nocite{Diop2017} Diop and Kengne (2017), \nocite{Lee2018asymptotic} Lee {\it et al.} (2018). Several authors have pointed out some restrictions of these procedures and proposed a Wald-type statistic based on a quasi likelihood estimators; see, among others papers \nocite{Lee2008test} Lee and Song (2008), \nocite{Kengne2012} Kengne (2012), \nocite{Diop2020b} Diop and Kengne (2020). Other procedures have been developed with the (conditional) least-squares estimator (see for instance  \nocite{Lee2005testa} Lee and Na (2005a), \nocite{Kang2009parameter} Kang and Lee (2009)) or the density power divergence estimator (see among others \nocite{Lee2005testb} Lee and Na (2005b),  \nocite{Kang2015robust} Kang and Song (2015)). 
 \nocite{Lee2003cusum} Lee {\it et al.} (2003) proposed a procedure for change-point detection in a large class of time series models. This procedure does not take into account the change-point alternative and does not ensure the consistency in power.

\medskip 

In this new contribution, we consider a multivariate continuous or integer valued process and deal with a general contrast, where the likelihood, quasi likelihood, conditional least-squares or density power divergence can be seen as a specific case.  
 %
%
 Let $\widehat{C} \big((Y_t)_{t \in T},\theta \big)$ be a contrast function defined for any segment $T \subset \{1,\cdots,n\}$ and $\theta \in \Theta$ by:
 \begin{equation}\label{C_def}
 \widehat{C} \big((Y_t)_{t \in T},\theta \big) = \sum_{t \in T} \widehat{\varphi}_t(\theta), 
  \end{equation} 
 where $\widehat{\varphi}_t$ depends on $Y_1,\cdots,Y_t$, and is such that, the minimum contrast estimator (MCE), computed on a segment $T \subset \{1,\cdots,n\}$ is given by 
 \begin{equation}\label{mce}
  \widehat{\theta}(T) \coloneqq  \underset{\theta\in \Theta}{\text{argmin}} \Big( \widehat{C} \big((Y_t)_{t \in T},\theta \big)  \Big).
  \end{equation}
In the sequel, we use the notation $\widehat{C} (T,\theta ) = \widehat{C} \big((Y_t)_{t \in T},\theta \big)$; and address the following issues.
\begin{enumerate}
\item[(i)] We propose a Wald-type statistic based on the MCE for testing H$_0$ against H$_1$. The asymptotic studies under the null and the alternative hypotheses show that, the test has correct size asymptotically and is consistent in power.
\item[(ii)] Application to a large class of multivariate causal processes is carried out. We provide sufficient conditions under which, the asymptotic results of the change-point detection hold. 
\item[(iii)] A general class of multivariate integer valued models is considered.  In the case where the conditional distribution belongs to the $m$-parameter exponential family, we provide sufficient conditions that ensure the existence of a stationary and ergodic $\tau$-weakly
dependent solution. The inference is carried out, and the consistency and the asymptotic normality of the Poisson quasi maximum likelihood estimator (PQMLE) are established. This inference question has been addressed by \nocite{Ahmad2016contribution} Ahmad (2016) with the equation-by-equation PQMLE, Lee  \textit{et al.} (2018)  for bivariate Poisson INGARCH model, \nocite{Cui2019flexible} Cui {\it et al.} (2019) for flexible bivariate Poisson integer-valued GARCH model, \nocite{Fokianos2020multivariate} Fokianos {\it et al.} (2020) for linear and log-linear multivariate Poisson autoregressive models. The model considered in Section \ref{Sec_class_MOD} appears to be more general and the conditions imposed for asymptotic studies  seem to be more straightforward. Also, we show that the asymptotic results of the change-point detection hold for this class of models.
\end{enumerate}

\medskip

The paper is structured as follows. Section 2 contains the general assumptions and the construction of the test statistic for change-point detection, as well as the main asymptotic results under H$_0$ and H$_1$. 
Section 3 is devoted to the application of the proposed change-point detection procedure to a general class of continuous valued processes. 
 Section 4 focuses on a general class of observation-driven integer-valued time series. Some simulation results are displayed in Section 5. 
Section 6 contains the proofs of the main results.


  \section{General change-point detection procedure}\label{Sec_detection}
  \subsection{Assumptions}
\noindent
 Throughout the sequel, the following norms will be used:
{\em
\begin{itemize}
 \item $ \|x \| \coloneqq  \sum_{i=1}^{p} |x_i| $ for any $x \in \mathbb{R}^{p}$ (with $p \in \N$);
%
\item $ \|x \| \coloneqq \underset{1\leq j \leq q}{\max} \sum_{i=1}^{p} |x_{i,j}| $ for any matrix $x=(x_{i,j}) \in M_{p,q}(\R)$; where $M_{p,q}(\R)$ denotes the set of matrices of dimension $p\times q$ with coefficients in $\R$;
\item  $\left\|g\right\|_{\mathcal K} \coloneqq \sup_{\theta \in \mathcal K}\left(\left\|g(\theta)\right\|\right)$ for any  compact set $\mathcal K \subseteq  \Theta$ and 
function $g:\mathcal K \longrightarrow   M_{p,q}(\R)$;
%
 %
\item $\left\|Y\right\|_r \coloneqq \E\left(\left\|Y\right\|^r\right)^{1/r}$ for any random vector  $Y$ with finite $r-$order moments. 
\end{itemize}
}  
  
\medskip

  Let $Y=\{Y_{t},\,t\in \Z \}$ be a multivariate continuous or integer valued process depending to a parameter $\theta^* \in \Theta$ and denote by $\mathcal{F}_{t-1}=\sigma\left\{Y_{t-1},\cdots \right\}$ the $\sigma$-field generated by the whole past at time $t-1$.
In the sequel, we assume that $(j_n)_{n\geq1}$ and $(k_n)_{n\geq1}$ are two integer valued sequences such that $j_n \leq k_n$, $k_n \rightarrow  \infty$ and $k_n - j_n \rightarrow  \infty $ as $n \rightarrow \infty$, and use the notation $T_{\ell,\ell'}=\{\ell,\ell+1,\cdots,\ell'\}$ for any $(\ell,\ell') \in \N^2$ such as $\ell \leq \ell'$.
We consider a segment $T_{j_n,k_n}$ and set the following assumptions for $(Y,\theta^*)$  under H$_0$.
  \begin{enumerate} 
     \item [(\textbf{A1}):] The process $Y=\{Y_{t},\,t\in \Z \}$ is  stationary and ergodic.
    \item [(\textbf{A2}):]  The minimum contrast estimator $\widehat{\theta}(T_{j_n,k_n})$ (defined in (\ref{mce})) converges  $a.s.$ to $\theta^*$.  
    \item [(\textbf{A3}):] For all $ t \in T_{j_n,k_n}$, the function $\theta  \mapsto \widehat{\varphi}_t(\theta)$ (see (\ref{C_def})) is continuously differentiable on $\Theta$, in addition,
     there exists a sequence of random function $(\varphi_t(\cdot))_{t \in \Z}$ such that, the mapping $\theta  \mapsto \varphi_t(\theta)$ is continuously differentiable on $\Theta$ and for all $\theta \in \Theta$, the sequence $(\partial \varphi_t(\theta)/ \partial \theta)_{t \in \Z}$ is  stationary and ergodic, satisfying:
  %
  %
   \begin{equation}\label{eq_assump_app_dif}
     \dfrac{1}{\sqrt{k_n - j_n}}  \sum_{t \in T_{j_n,k_n}} \Big \| \dfrac{\partial}{\partial \theta} \widehat{\varphi}_t(\theta) 
     - \dfrac{\partial}{\partial \theta} \varphi_t(\theta) \Big \|_\Theta 
     = o_P(1) ~ \text{ and } ~
     \E \Big\| \dfrac{\partial}{\partial \theta} \varphi_t(\theta) \Big \|^2_\Theta < \infty.  
  \end{equation} 
  Furthermore, $\left(\frac{\partial}{\partial \theta} \varphi_t(\theta^*),\mathcal{F}_{t}\right)_{t \in \mathbb{Z}}$ is a stationary ergodic, square integrable martingale difference sequence with covariance $G  =  \E \big[ \dfrac{\partial \varphi_0(\theta^*)}{\partial \theta}  \dfrac{\partial \varphi_0(\theta^*)}{\partial \theta^T}\big] $ assumed to be positive definite.
 
  \item [(\textbf{A4}):]  For all $ t \in T_{j_n,k_n}$, the function $\theta  \mapsto \widehat{\varphi}_t(\theta)$  is 2 times continuously differentiable on $\Theta$, moreover, under the assumption (\textbf{A3}), the function $\theta  \mapsto \frac{\partial \varphi_t(\theta)}{\partial \theta}$ is continuously differentiable on $\Theta$, such that, 
the sequence $(\partial^2 \varphi_t(\theta)/ \partial \theta \partial \theta^T)_{t \in \Z}$ is  stationary and ergodic, satisfying:  
  %
   %
 \begin{equation}\label{eq_assump_app_dif2}
 \E \Big\|  \dfrac{\partial^2 \varphi_t(\theta)}{\partial \theta \partial \theta^T} \Big \|_\Theta < \infty, ~ ~ 
    \Big \| \dfrac{1}{ k_n - j_n}  \sum_{t \in T_{j_n,k_n}}  \dfrac{\partial^2 \widehat{\varphi}_t(\theta)}{\partial \theta \partial \theta^T} - \E \Big(  \dfrac{\partial^2 \varphi_0(\theta)}{\partial \theta \partial \theta^T} \Big)  \Big \|_\Theta 
     = o(1),
  \end{equation} 
  %
  %
  and the matrix $F =  \E \big[ \dfrac{\partial^2 \varphi_0(\theta^*)}{\partial \theta \partial \theta^T}  \big]$ is invertible. In addition, $\widehat{\theta}(T_{j_n,k_n})$ is asymptotically normal; that is,
 \begin{equation}\label{Assymp_normal}
  \sqrt{k_n - j_n}\left(\widehat{\theta}(T_{j_n,k_n})-\theta^*\right) \limiteloin \mathcal{N}(0,\Omega)~ \text{ with } ~\Omega \coloneqq F ^{-1} G  F^{-1} . 
  \end{equation}
  \end{enumerate}

  \medskip
  \noindent
  For any $\ell, \ell' \in \N$ with $\ell \leq \ell'$, define  the  matrices:
\[
  \widehat G(T_{\ell,\ell'})  = 
  \frac{1}{\ell'-\ell+1}\sum_{t \in T_{\ell,\ell'}}  \dfrac{\partial}{\partial \theta} \widehat{\varphi}_t(\widehat{\theta}(T_{\ell,\ell'})) \dfrac{\partial}{\partial \theta^T} \widehat{\varphi}_t(\widehat{\theta}(T_{\ell,\ell'}))
 ~~\text{and}~~   
   \widehat F(T_{\ell,\ell'})  = 
    \frac{1}{\ell'-\ell+1}\sum_{t \in T_{\ell,\ell'}} \dfrac{\partial^2}{\partial \theta \partial \theta^T} \widehat{\varphi}_t(\widehat{\theta}(T_{\ell,\ell'})) .
\]
 According to (\textbf{A1})-(\textbf{A4}),  $\widehat F(T_{j_n,k_n})$ and $\widehat G(T_{j_n,k_n})$ converges almost surely to $F $ and $G$, respectively. 
Therefore, $\widehat F(T_{j_n,k_n})^{-1} \widehat G(T_{j_n,k_n})   \widehat F(T_{j_n,k_n})^{-1}$ is a consistent estimator of the covariance matrix $\Omega $. 

 \subsection{Change-point test and asymptotic results} 
 We derive a retrospective test procedure based on the MCE of the parameter. 
For all $n \geq 1$, define the matrix $\widehat{\Omega}(u_n)$ and  the subset $\mathcal T_n$ by  
\begin{equation*}
\widehat{\Omega}(u_n)=\frac{1}{2} 
\big[
\widehat F(T_{1,u_n})  \widehat G(T_{1,u_n})^{-1}   \widehat F(T_{1,u_n}) +
\widehat F(T_{u_n+1,n})  \widehat G(T_{u_n+1,n})^{-1}  \widehat F(T_{u_n+1,n})
\big]
~\text{ and }~ \mathcal T_n= [v_n , n-v_n]\cap \N, 
\end{equation*}
where $(u_n,v_n)_{n\geq 1}$ is a bivariate integer valued sequence such that:  $ (u_n,v_n) =o(n)$ and $ u_n,v_n \limiten +\infty$.\\
 For any $1< k < n$, let us introduce 
 \[ \widehat{Q}_{n,k}=\frac{\left(k(n-k)\right)^{2}}{n^{3}}
 \big(\widehat{\theta}(T_{1,k})-\widehat{\theta}(T_{k+1,n})\big)^T
 \widehat{\Omega}(u_n)
  \big(\widehat{\theta}(T_{1,k})-\widehat{\theta}(T_{k+1,n})\big).
  \]
Therefore, consider the following test statistic:
\begin{equation}\label{Def_stat_test}
  \widehat{Q}_n=\max_{k \in \mathcal T_n}\big(\widehat{Q}_{n,k}\big).
\end{equation}
The construction of this statistic follows the approach of Doukhan and Kengne (2015); 
that is, $\widehat{Q}_{n,k}$ evaluates a distance between $\widehat{\theta}(T_{1,k})$ and $\widehat{\theta}(T_{k+1,n})$ for all $k \in \mathcal  T_n$.
The null hypothesis  H$_0$ will thus be rejected if there exists a time $k \in \mathcal  T_n$ such that this distance is too large.

 \medskip
\noindent
The following theorem gives the asymptotic behavior of the test statistic under H$_0$. 
\begin{thm}\label{th1_test}
Under H$_0$ with $\theta^* \in \overset{\circ}{\Theta}$, assume that (\textbf{A1})-(\textbf{A4}) hold for $(Y,\theta^*)$.
Then,
\begin{equation}\label{res_th1_test}
\widehat{Q}_n \limiteloin \sup_{0\leq  \tau \leq 1}\left\|W_d(\tau)\right\|^{2}, 
\end{equation}
where $W_d$ is a $d$-dimensional Brownian bridge.
\end{thm}
For a nominal level $\alpha \in (0,1)$, the critical region of the test is then $(\widehat{Q}_{n}>c_{d,\alpha})$, where $c_{d,\alpha}$ is the $(1-\alpha)$-quantile of the distribution
of $\underset{0\leq  \tau \leq 1}{\sup}\left\|W_d(\tau)\right\|^{2}$. The critical values  $c_{d,\alpha}$ can be easily obtained through a Monte Carlo simulation; see for instance, Lee {\it et al.} (2003).
%

  \medskip
  
    \medskip
 \noindent Under the alternative hypothesis, we consider the following additional condition for the break instant.

  \medskip
\noindent {\bf Assumption B}: {\em There exists $\tau^* \in (0,1)$ such that $t^*=[n\tau^*]$,  where $[x]$ denotes the integer part of $x$.}

\medskip

\medskip
\noindent
 We obtain the following main result under H$_1$.
\begin{thm}\label{th2_test}
Under H$_1$ with $\theta^*_1$ and $\theta^*_2$ belonging to  $\overset{\circ}{\Theta}$, assume that  (\textbf{A1})-(\textbf{A4}) hold for $(Y^{(1)},\theta^*_1)$ and $(Y^{(2)},\theta^*_2)$.  If Assumption {\bf B} is satisfied, then

\begin{equation}\label{res_th2_test}
 \widehat{Q}_n  \limiteproban +\infty .
 \end{equation}
\end{thm}

 \medskip
 \noindent 
 In the next two sections, we will  detail  some examples of classes of multivariate time series with a quasi likelihood contrast function. 
  We also show that, under some regularity conditions, the  general assumptions required for Theorems \ref{th1_test} and \ref{th2_test} are satisfied for these classes. 
 Let us stress that, the scope of the proposed procedure is quite extensive and is not only restricted to the examples below. This procedure can be applied for instance, for change-point detection in models with  exogenous covariates (see \nocite{Diop2021inference} Diop and Kengne (2021), \nocite{Aknouche2021count} Aknouche and Francq (2021)), for integer valued time series with negative binomial quasi likelihood contrast (see \nocite{Aknouche2018} Aknouche {\it et al.} (2018)) or with density power divergence contrast (see \nocite{Kim2020robust} Kim and Lee (2020)), for general time series model with the conditional least-squares contrast (see \nocite{Klimko1978conditional} Klimko and Nelson (1978)). In fact, one can easily see in these papers that, the assumptions (\textbf{A1})-(\textbf{A4}) hold.

 
 \section{Application to a class of multidimensional causal processes}\label{Sec_class_AC}
Let $\{Y_{t},\,t\in \Z \}$ be a multivariate  time series of dimension $m \in \N$. 
 For any  $\mathcal T \subseteq \Z$ and $\theta \in \Theta$, consider the general class of causal processes defined by 
  \medskip
              
   \textbf{Class} $\mathcal{AC}_{\mathcal T}(M_{\theta},f_{\theta})$: A process $\{Y_{t},\,t\in \mathcal T \}$ belongs to $\mathcal{AC}_{\mathcal T}(M_\theta,f_\theta)$ if it satisfies:
   \begin{equation}\label{Model_AC_mf} 
     Y_t =M_{\theta}(Y_{t-1}, Y_{t-2}, \ldots) \cdot \xi_t + 
      f_{\theta}(Y_{t-1}, Y_{t-2}, \ldots)~~\forall t \in \mathcal T,
   \end{equation}
  where 
   $M_{\theta}(Y_{t-1}, Y_{t-2}, \ldots)$ is a $m \times p$ 
  random matrix having almost everywhere ($a.e$) full rank $m$, $f_{\theta}(Y_{t-1}, Y_{t-2}, \ldots)$ is a $\R^m$-random vector 
  and 
   $(\xi_t)_{t \in \Z}$ is a $\R^p$-random vector independent, identically distributed (\textit{i.i.d}) satisfying $\xi_t=(\xi^{(k)}_t)_{1\leq k \leq p}$ with $\E\big[\xi^{(k)}_0 \xi^{(k')}_0\big]=0$ for $k \neq k'$ and $\E\big[\xi^{{(k)}^2}_0\big]=\text{Var}(\xi^{(k)}_0)=1$ for $1\leq k \leq p$.
   $M_\theta(\cdot)$ and $f_\theta(\cdot)$ are assumed to be known up to the parameter $\theta$. 
    This class  has  been studied in \nocite{Doukhan2008} Doukhan and Wintenberger (2008), Bardet and Wintenberger (2009).
    
    \medskip   
    \noindent
Now, assume that $(Y_1,\ldots,Y_n)$ is a trajectory generated from one or two processes satisfying  (\ref{Model_AC_mf}). We would like to carry out the change-point test presented in Section \ref{sec_intro}.     
For all $t \in \Z$, denote by $\mathcal{F}_{t}=\sigma(Y_s,\, s \leq t)$ the $\sigma$-field generated by the whole past at time $t$. 
For any segment $T \subset \{1,\cdots,n\}$ and $\theta \in \Theta$, we define the contrast function based on the conditional Gaussian quasi-log-likelihood given by (up to an additional constant)
 \begin{equation}\label{def_C_T.theta_AC} 
 \widehat  C(T,\theta)=  \frac{1}{2}\sum_{t \in T} \widehat \varphi_t(\theta)
 ~~~ \text{with}~~~
        \widehat \varphi_t(\theta) = (Y_t-\widehat f^t_\theta)^T(\widehat H^t_\theta)^{-1}(Y_t-\widehat f^t_\theta)+\log( \det (\widehat H^t_\theta)),
\end{equation}      
where $\widehat f^t_\theta := f_\theta(Y_{t-1},\ldots,Y_{1},0,\cdots)$, 
$\widehat M^t_\theta := M_\theta(Y_{t-1},\ldots,Y_{1},0,\cdots)$, $\widehat H^t_\theta := \widehat M^t_\theta  (\widehat{M}^{t}_\theta)^T$.
 %
 Thus, the MCE  computed on $T$ is defined by
 \begin{equation}\label{mce_AC} 
 \widehat{\theta}(T)= \underset{\theta\in T}{\text{argmin}} \big(\widehat C(T,\theta) \big).
\end{equation}

\medskip
\noindent 
Let $\Psi_\theta$ be a generic symbol for any of the functions $f_\theta$, $M_\theta$ or $H_\theta=M_\theta M_\theta^T$ and $\mathcal K \subseteq \Theta$  be a compact subset.  
  To study the stability properties of the class (\ref{Model_AC_mf}),  Bardet and Wintenberger (2009) imposed the following classical Lipschitz-type conditions on the function $\Psi_\theta$.
  
\medskip
   
    \noindent 
    \textbf{Assumption} \textbf{A}$_i (\Psi_\theta,\mathcal K)$ ($i=0,1,2$):
    For any $y \in (\R^m)^{\infty}$, the function $\theta \mapsto \Psi_\theta(y)$ is $i$ times continuously differentiable on $\mathcal K$  with $ \big\| \frac{\partial^i \Psi_\theta(0)}{\partial \theta^i}\big\|_{\mathcal K}<\infty $; 
    and
      there exists a sequence of non-negative real numbers $(\alpha^{(i)}_{k}(\Psi_\theta,\mathcal K))_{k \in \N} $  satisfying:
     $ \sum\limits_{k=1}^{\infty} \alpha^{(i)}_{k}(\Psi_\theta,\mathcal K) <\infty$, for $i=0, 1, 2$;
   such that for any  $x, y \in (\R^m)^{\infty}$,
  \[  \Big \| \frac{\partial^i}{ \partial \theta^i} \Psi_\theta(x)-\frac{\partial^i}{\partial\theta^i} \Psi_\theta(y) \Big \|_{\mathcal K}
  \leq  \sum\limits_{k=1}^{\infty}\alpha^{(i)}_{k}(\Psi_\theta,\mathcal K) \|x_k-y_k\|, 
  \]
where 
$x,y,x_k,y_k$ are respectively replaced by $xx^T$, $yy^T$, $x_kx^T_k$, $y_ky^T_k$ if $\Psi_\theta=H_\theta$. 
\medskip

\noindent 
    For $r\geq 1$, define the set 
\begin{multline*}\label{Set_Theta(r)}
\Theta(r) = \big\{
 \theta \in \R^d \, \big / \, \textbf{A}_0 (f_\theta,\{\theta\}) \   \text{and}\ \textbf{A}_0 (M_\theta,\{\theta\})   \    \text{hold with} 
 ~ \sum\limits_{k=1}^{\infty} \left\{\alpha^{(0)}_{k}(f_\theta,\{\theta\}) + \|\xi_0\|_r \alpha^{(0)}_{k}(M_\theta,\{\theta\}\right\} 
 <1
\big\}\\
\bigcup 
\big\{
 \theta  \in \R^d \ \big / \ f_\theta=0 \text{ and } \textbf{A}_0 (H_\theta,\{\theta\})  \text{ holds with } 
 \|\xi_0\|^2_r  \sum\limits_{k=1}^{\infty} \alpha^{(0)}_{k}(H_\theta,\{\theta\}) 
 <1
\big\}.
\end{multline*}

\medskip

\noindent 
The following regularity conditions are also considered in Bardet and Wintenberger (2009) to assure the consistency and the asymptotic normality of $\widehat{\theta}(T_{1,n})$ under H$_0$.

\medskip
\noindent
  ($\mathcal{AC}.\textbf{A0}$): For all  $\theta \in \Theta$ and some $t \in \Z$, 
 $ \big( f^t_{\theta^*}= f^t_{\theta}  \ \text{and} \ H^t_{\theta^*}= H^t_{\theta} \ \ a.s. \big) \Rightarrow ~ \theta= \theta^*$.
 
\medskip
\noindent
  ($\mathcal{AC}.\textbf{A1}$): $\exists  \underline{H}>0$ such that $\displaystyle \inf_{ \theta \in \Theta} \det \left(H_\theta(y)\right)  \geq \underline{H}$,  for all $y \in (\R^m)^{\infty}$.
  
  \medskip
  \noindent
  ($\mathcal{AC}.\textbf{A2}$):  
 $
 \alpha^{(i)}_{k}(f_\theta,\Theta)+\alpha^{(i)}_{k}(M_\theta,\Theta)
 +\alpha^{(i)}_{k}(H_\theta,\Theta) = O(k^{-\gamma})
$
 for $i=0,1,2$ and some $\gamma >3/2$.
 
  \medskip
  \noindent
  ($\mathcal{AC}.\textbf{A3}$): One of the families $\big(\frac{\partial f^0_{\theta^*}}{\partial \theta_i}\big)_{1\leq i \leq d}$ or 
    $\big(\frac{\partial H^0_{\theta^*}}{\partial \theta_i}\big)_{1\leq i \leq d}$ is  $a.e$ linearly
independent.\\

\medskip
\noindent
Under \textbf{A}$_0(\Psi_\theta,\Theta)$ (for $\Psi_\theta=f_\theta, M_\theta, H_\theta$) with $\theta^* \in \Theta \cap \Theta(1)$, \nocite{Bardet2009} Bardet and Wintenberger (2009) established the existence of a $\tau$-weakly dependent strictly stationary and ergodic solution to   the class $\mathcal{AC}_{\Z}(M_{\theta^*},f_{\theta^*})$;  which shows that the assumption (\textbf{A1}) holds.
Under H$_0$, if  \textbf{A}$_0(f_\theta,\Theta)$, \textbf{A}$_0(M_\theta,\Theta)$ (or \textbf{A}$_0(M_\theta,\Theta)$) and ($\mathcal{AC}.\textbf{A0}$)-($\mathcal{AC}.\textbf{A2}$) hold with $\theta^* \in \Theta \cap \Theta(2)$, then $\widehat{\theta}(T_{j_n,k_n}) \limitepsn \theta^*$ (from Theorem 1 of Bardet and Wintenberger (2009)).
  Therefore, (\textbf{A2}) is satisfied. 

\medskip
\noindent
Let us define   
 \begin{equation}\label{def_fi.t_AC}
\varphi_t(\theta):= (Y_t- f^t_\theta)^T(H^t_\theta)^{-1}(Y_t- f^t_\theta)+\log( \det (H^t_\theta))
 \end{equation}
  with $f^t_\theta := f_\theta(Y_{t-1},\ldots)$,   
  $M^t_\theta := M_\theta(Y_{t-1},\ldots)$ and $H^t_\theta := M^t_\theta ({M^t_\theta})^T$.\\
 Under H$_0$, \textbf{A}$_i(f_\theta,\Theta)$, \textbf{A}$_i(M_\theta,\Theta)$ (or \textbf{A}$_i(M_\theta,\Theta)$) for $i=0,1,2$ and ($\mathcal{AC}.\textbf{A0}$)-($\mathcal{AC}.\textbf{A3}$) with $\theta^* \in \overset{\circ}{\Theta} \cap \Theta(4)$,  Bardet and Wintenberger (2009) have proved that $\widehat{\theta}(T_{j_n,k_n})$ is asymptotically normal. 
 Then, using the sequence of functions $\left(\varphi_t(\cdot)\right)_{t \in \Z}$ defined in (\ref{def_fi.t_AC}), one can see that 
the assumptions (\textbf{A3}) and (\textbf{A4}) also hold.   
For the condition (\ref{eq_assump_app_dif}) and those imposed on the sequence  $(\frac{\partial}{\partial \theta} \varphi_t(\theta^*),\mathcal{F}_{t})_{t \in \mathbb{Z}}$ in (\textbf{A3}), see the proof of their Theorem 2.
 Thus, under the null hypothesis,  all the  required assumptions  (\textbf{A1})-(\textbf{A4}) are verified for $(Y,\theta^*)$; which assures that Theorem \ref{th1_test} applies to this class of models. 
Note that, by the same arguments, one can also see that, these assumptions hold for $(Y^{(1)},\theta^*_1)$ and $(Y^{(2)},\theta^*_2)$ under H$_1$. Therefore, Theorem  \ref{th2_test} also applies to this class.

%

\section{Inference and application in general multivariate count process}\label{Sec_class_MOD}

\subsection{Model formulation and inference}
 Consider a multivariate count time series $\{Y_{t}= (Y_{t,1},\cdots,Y_{t,m})^T,\,t\in \Z \}$ with value in $\N_0^{m}$ (with $m \in \N$, $\N_0 = \N \cup \{ 0 \}$) and denote by $\mathcal{F}_{t-1}=\sigma\left\{Y_{t-1},\cdots  \right\}$ the $\sigma$-field generated by the whole past at time $t-1$. 
 For any $\mathcal T \subseteq \Z$ and $\theta \in \Theta$, define the class of multivariate observation-driven integer-valued time series given by
 
 \medskip
 {\bf Class} $\mathcal{MOD}_{\mathcal T }(f_{\theta})$: The multivariate count process $Y=\{Y_{t},\,t\in \mathcal T  \}$ belongs to $\mathcal{MOD}_{\mathcal T }(f_{\theta})$ if it satisfies:
  \begin{equation} \label{Model_MOD}
             \E(Y_t|\mathcal{F}_{t-1})=f_{\theta}(Y_{t-1},Y_{t-2},\cdots) ~~ \forall t \in \mathcal{T} ,
   \end{equation}
where $f_{\theta}(\cdot)$ is a measurable multivariate function with  non-negative components, assumed to be known up to the  parameter $\theta$.

 \medskip
\noindent
In this section, it is assumed that any $\{Y_{t} ,\,t\in \Z\}$ belonging to $\mathcal{MOD}_{\mathcal T }(f_{\theta})$ is a stationary and ergodic process (i.e, the condition (\textbf{A1}) imposed for the change-point detection holds) satisfying: 
 \begin{equation}\label{moment}
    \exists C>0, \epsilon >0, \text{ such that } \forall t \in \Z, ~ ~ \|Y_{t}\|_{1+\epsilon} <C. 
   \end{equation}
   
   \medskip  
 \noindent 
Let $(Y_{1},\ldots,Y_{n})$ be observations generated from $\mathcal{MOD}_{\Z}(f_{\theta^{*}})$ with $\theta^* \in \Theta$. 
The conditional Poisson quasi log-likelihood computed on $\{1,\cdots,n\}$ is given  by (up to a constant)
 \[
  L_n(\theta) \coloneqq 
   \sum_{t=1}^{n} \ell_t(\theta) ~ \text{ with }~
  \ell_t(\theta) = \sum_{i=1}^{m} \left(Y_{t,i}\log \lambda_{t,i}(\theta)- \lambda_{t,i}(\theta)\right),
  \]
  where $ \lambda_t(\theta) \coloneqq \big(\lambda_{t,1}(\theta), \cdots, \lambda_{t,m}(\theta) \big) =f_\theta(Y_{t-1}, Y_{t-2}, \cdots )$.
 An approximated conditional quasi log-likelihood is given by
\begin{equation*}
\widehat{L}_n( \theta) \coloneqq  \sum_{t=1}^{n} \widehat{\ell}_t(\theta) ~ \text{ with }~
 \widehat{\ell}_t(\theta) =   \sum_{i=1}^{m} \left(Y_{t,i} \log \widehat{\lambda}_{t,i}(\theta)- \widehat{\lambda}_{t,i}(\theta)\right),
 \end{equation*}
  where $ \widehat{\lambda}_t(\theta) \coloneqq \big(\widehat{\lambda}_{t,1}(\theta), \cdots, \widehat{\lambda}_{t,m}(\theta) \big)^T =  f_\theta(Y_{t-1}, \cdots, Y_{1},0,\cdots)$.
 Therefore, the Poisson quasi-maximum likelihood estimator (QMLE) of $ \theta^*$  is defined by
 \begin{equation*}
  \widehat{\theta}_n \coloneqq  \underset{\theta\in \Theta}{\text{argmax}} \big(\widehat{L}_n(\theta)\big).
  \end{equation*}
 Note that, under the assumption of independence among components and conditionally Poisson distributed, this Poisson QMLE is equivalent to the maximum likelihood estimator. Let us highlight that, we deal with an arbitrary dependence among components and arbitrary conditional distribution; that is, the distribution of the components could differ from each other.
  \medskip
  
  \noindent
  For a process $\{Y_t,\, t \in \Z\}$ belonging to $\mathcal{MOD}_{\Z}(f_{\theta^{*}})$, we set the following assumptions in order to establish the consistency and the asymptotic normality of the Poisson QMLE. 

   \medskip
   \noindent
  \textbf{Assumption} \textbf{A}$_i (\Theta)$ ($i=0,1,2$):
    For any $y \in \big(\mathbb{N}_0^{m} \big)^{\infty}$, the function $\theta \mapsto f_\theta(y)$ is $i$ times continuously differentiable on $\Theta$  with $ \left\| \partial^i f_\theta(0)/ \partial \theta^i\right\|_\Theta<\infty $; 
    and
      there exists a sequence of non-negative real numbers $(\alpha^{(i)}_k)_{k\geq 1} $ satisfying
     $ \sum\limits_{k=1}^{\infty} \alpha^{(0)}_k <1 $ (or $ \sum\limits_{k=1}^{\infty} \alpha^{(i)}_k <\infty $ for $i=1, 2$);
   such that for any  $y, y' \in \big(\mathbb{N}_0^{m} \big)^{\infty}$,
  \[ 
   \Big \| \frac{\partial^i f_\theta(y)}{ \partial \theta^i}-\frac{\partial^i f_\theta(y')}{\partial\theta^i} \Big \|_\Theta
  \leq  \sum\limits_{k=1}^{\infty}\alpha^{(i)}_k \|y_k-y'_k\|.
  \]

  \medskip
  \noindent
 ($\mathcal{MOD}.\textbf{A0}$): 
     For all  $\theta \in \Theta$,
 $ \big( f_{\theta^*}(Y_{t-1}, Y_{t-2}, \cdots) \equalpsn  f_{\theta}(Y_{t-1}, Y_{t-2}, \cdots)  ~ \text{ for some } t \in \Z \big) \Rightarrow ~ \theta^* = \theta$; 
 moreover, $\exists  \underline{c}>0$ such that $ f_\theta(y)  \geq \underline{c} \textbf{1}_m$ componentwise, for all $\theta \in \Theta$, $ y \in  \big( \N_0^{m} \big)^{\infty} $, where $\textbf{1}^T_m=(1,\cdots,1)$ is a vector of dimension $m$.

 \medskip
 \noindent 
 ($\mathcal{MOD}.\textbf{A1}$):  $\theta^* $ is an interior point of $\Theta \subset \mathbb{R}^{d}$.

 \medskip
 \noindent
  ($\mathcal{MOD}.\textbf{A2}$):  The family $ \big(\frac{\partial \lambda_{t} (\theta^* )}{\partial    \theta_i} \big)_{1\leq i \leq d} $  is $a.e.$  linearly independent.

 \medskip
 
  \medskip
   \noindent
 Proposition \ref{prop1} below establishes the existence of a stationary and ergodic solution of the model (\ref{Model_MOD}) for the $m$-parameter exponential family conditional distribution.
 Consider a $m$-dimensional process $\{Y_t, ~ t \in \Z \}$ satisfying 
 \begin{equation}\label{Model_exp_f}
  Y_t|\mathcal{F}_{t-1} \sim p(y|\eta_t) ~ \text{ with } ~ \lambda_t(\theta) := \E(Y_t|\mathcal{F}_{t-1})=f_{\theta}(Y_{t-1},Y_{t-2},\cdots)  
 \end{equation}
 where $p(\cdot | \cdot )$ is a multivariate discrete distribution belonging to the $m$-parameter exponential family; that is
 \[  p(y|\eta) = \exp\{\eta^T y - A(\eta) \} h(y), ~ y \in \N^m_0 \] 
 where $\eta$ is the natural parameter (i.e. $\eta_t$ is the natural parameter of the distribution of $Y_t|\mathcal{F}_{t-1}$) and $A(\eta)$, $h(y)$ are known functions. It is assumed that the function $\eta \mapsto A(\eta)$ is twice continuously differentiable on the natural parameter space; therefore, the mean and variance of this distribution are $\partial A(\eta)/ \partial \eta$ and  $\partial^2 A(\eta)/ \partial \eta^2$, respectively. See \nocite{Khatri1983}  Khatri (1983) for more details on such class of distribution.
 For the model (\ref{Model_exp_f}), it holds that 
 \[ \E(Y_t|\mathcal{F}_{t-1})=f_{\theta}(Y_{t-1},Y_{t-2},\cdots) = \dfrac{\partial A(\eta_t)}{ \partial \eta} .\]
\begin{prop}\label{prop1}
Assume that \textbf{A}$_0(\Theta)$  holds. Then,  there exists a $\tau-weakly$ dependent,  stationary and ergodic solution $\{ Y_t, ~ t \in \Z \}$ to (\ref{Model_exp_f}), satisfying $\E\|Y_t \| < \infty$.
\end{prop}

 \medskip
   \noindent
In the sequel, we deal with the more general class of model (\ref{Model_MOD}), where the distribution of $Y_t|\mathcal{F}_{t-1}$ may be outside the $m$-parameter exponential family.
 The following theorem shows that the Poisson QMLE for the class of  models (\ref{Model_MOD}) is strongly consistent.
\begin{thm}\label{th1}
Assume that \textbf{A}$_0(\Theta)$, ($\mathcal{MOD}.\textbf{A0}$)  and (\ref{moment}) (with $\epsilon \geq  1$) hold with 
\begin{equation}\label{cond_th1}
 \alpha^{(0)}_{k}= \mathcal{O}(k^{-\gamma}) \text{ for some } \gamma >3/2.
\end{equation}
Then  
\begin{equation*}\label{Cons_theta}
\widehat{\theta}_{n} \limitepsn \theta^*.
\end{equation*}
\end{thm}

 \medskip
  \noindent
 For any $t \in \Z$ and $\theta \in \Theta$, denote $\Gamma_t (\theta):=(Y_t -\lambda_{t}(\theta))(Y_t -\lambda_{t}(\theta))^T$  and $D_t(\theta)$ the $m \times m$ diagonal matrix with the $i$th diagonal element is equal to $\lambda_{t,i}(\theta)$ for any $i=1,\cdots,m$.
 From the assumption ($\mathcal{MOD}.\textbf{A0}$), the matrix $D_t(\theta)$ is $a.s.$ positive definite. 
Combining all the regularity assumptions and notations given above, we obtain the asymptotic normality of the Poisson QMLE, as shown in
the following theorem.
\begin{thm}\label{th2}
Assume that \textbf{A}$_i(\Theta)$ ($i=0,1,2$), ($\mathcal{MOD}.\textbf{A0}$)-($\mathcal{MOD}.\textbf{A2}$) and (\ref{moment}) (with $\epsilon \geq 3$) hold with 
	 \begin{equation}\label{cond.th2}
  \alpha_k^{(0)} +\alpha_k^{(1)}= \mathcal O (k^{-\gamma}) ~~ for~ some~ \gamma>3/2,
 \end{equation}
 then 
 \[ \sqrt{n}(\widehat{\theta}_n-\theta^*) \limiteloin \mathcal{N}(0,\Sigma)~ \text{ with } ~\Sigma \coloneqq J^{-1}_{\theta^*} I_{\theta^*} J^{-1}_{\theta^*}, \]
 where
 \[
 J_{\theta^*} = \big[\frac{\partial \lambda^T_{0}(\theta^* )}{ \partial \theta} D^{-1}_0(\theta^*)   \frac{\partial \lambda_{0}(\theta^* )}{ \partial \theta^T } \big] 
           ~~and~~  \it
           I_{\theta^*}= \E \big[ \frac{\partial \lambda^T_{0}(\theta^* )}{ \partial \theta} D^{-1}_0(\theta^*) \Gamma_0 (\theta^*) D^{-1}_0(\theta^*) \frac{\partial \lambda_{0}(\theta^* )}{ \partial \theta^T} \big].
 \]
\end{thm}


\subsection{Change-point detection}
Now, assume that the trajectory $(Y_1,\ldots,Y_n)$ is  generated from one or two processes satisfying the general model (\ref{Model_MOD}) and consider the change-point test problem of Section \ref{sec_intro}.  
Let us define the contrast function based on the conditional Poisson quasi-log-likelihood for any segment $T \subset \{1,\cdots,n\}$ and $\theta \in \Theta$:
\begin{equation}\label{def_C_T.theta_MOD}
\widehat{C}(T, \theta) \coloneqq  \sum_{t \in T} \widehat{\varphi}_t(\theta) ~~ \text{with}~~
 \widehat{\varphi}_t(\theta) = -\widehat{\ell}_t(\theta) \text{ for all } t \in \Z. 
 \end{equation}
  %
  %
 Thus, the MCE computed on $T$ is given by 
 \begin{equation}\label{mce_MOD}
  \widehat{\theta}(T) \coloneqq  \underset{\theta\in \Theta}{\text{argmin}} \big(\widehat{C}(T, \theta)\big).
  \end{equation}

\medskip
\noindent 
Under the null hypothesis, the assumption (\textbf{A2}) holds from Theorem \ref{th1}. 
Letting $\varphi_t(\theta):=  -\ell_t(\theta)$ for all $ t \in \Z$ and $\theta \in \Theta$, one can see that, (\textbf{A3}) and (\textbf{A4}) are also  satisfied from Theorem \ref{th2}. The relation (\ref{eq_assump_app_dif}) in (\textbf{A3}) holds from Lemma \ref{lem2} (i) (see below) and the proof of Lemma \ref{lem3} ($a$), whereas the relation (\ref{eq_assump_app_dif2})  in (\textbf{A4}) holds from Lemme \ref{lem2} (ii), Lemma \ref{lem3} (c) and by applying the uniform law of large numbers. See also Lemma \ref{lem3}(b) for the required properties about  the sequence $(\frac{\partial}{\partial \theta} \varphi_t(\theta^*),\mathcal{F}_{t})_{t \in \mathbb{Z}}$. 
Hence,   in absence of change, all the conditions of Theorem \ref{th1_test} are verified for $(Y,\theta^*)$; which assures that the first result about the asymptotic behavior of the test statistic $\widehat{Q}_n$  applies  to the class of models (\ref{Model_MOD}).  
Under the change point alternative H$_1$, one can go along similar lines to verify that (\textbf{A1})-(\textbf{A4}) are satisfied for $(Y_1,\theta^*_1)$ and $(Y_2,\theta^*_2)$. This shows that Theorem \ref{th1_test} can also be applied to this class.
 



\section{Simulation study}
In this section, we investigate the performance (level and power) of the test statistic through two examples of two-dimensional processes, with sample size $n = 500, 1000$ and the nominal level $\alpha = 0.05$. 
For a sample size $n$, the statistic $\widehat{Q}_{n}$ will be computed  with $u_n=[\left(\log(n)\right)^{2}]$ and $v_n=[\left(\log(n)\right)^{5/2}]$. 
The procedure is implemented on the {\it R} software (developed by the CRAN project).

\medskip
\noindent
Let us consider the following models.
\begin{itemize}
	\item \emph{\bf A bivariate AR$(1)$ model.}\\ 
Consider the two-dimensional AR$(1)$ model expressed as
\begin{equation}\label{AC_SIM}
Y_t={A_0}Y_{t-1}+\xi_t ~
\text{ for all } t \in \Z,
\end{equation}
where $Y_{t}= (Y_{t,1},Y_{t,2})^T$, ${A_0}=(a_{i,j})_{i,j=1,2}$ is a $2 \times 2$ matrix with $\left\|{A_0}\right\|<1$ and $(\xi_t=(\xi_{t,1},\xi_{t,2})^T)_{t \in \Z}$ is a white noise  satisfying the conditions of the  general class (\ref{Model_AC_mf}).  
The parameter of the model is denoted by $\theta_0=(a_{1,1},a_{1,2},a_{2,1},a_{2,2})$. At the nominal level $\alpha = 0.05$, the critical value of the test is therefore $c_{4,\alpha} \approx 3.452$ (see  Lee \textit{et al.} (2003)). 
The performance will be evaluated in cases where the innovation $(\xi_{t})_{t \in \Z}$  is obtained from Student distributions with degrees of freedom  higher than $5$.
In the scenarios of change, we assume that the parameter  $\theta_0$
 (i.e, the matrix ${A_0}$) changes to $\theta_1$ which will be represented by a matrix denoted ${A_1}$.

\item \emph{\bf A bivariate INARCH$(1)$ model.} \\
Assume that $\{Y_{t}= (Y_{t,1},Y_{t,2})^T,\,t\in \Z \}$ is a count time series with value in $\N_0^{2}$, where $\{Y_{t,1},\,t\in \Z \}$ and $\{Y_{t,2},\,t\in \Z \}$ are two processes with conditional distribution following a Poisson distribution and a negative binomial distribution, respectively. More precisely,
\begin{equation}\label{MOD_SIM}
\left\{
\begin{array}{l}
  Y_{t,1} | \mathcal{F}_{t-1}  \sim \textrm{Poisson}(\lambda_{t,1})
  \\
\rule[0cm]{0cm}{.6cm}
Y_{t,2} | \mathcal{F}_{t-1}  \sim \textrm{NB}(r,r/(r+\lambda_{t,2}))
\end{array}
\right.
~~~~\text{ with } ~~\lambda_{t}:= 
(\lambda_{t,1},\lambda_{t,2})^T = {d_0}+{B_0}Y_{t-1},
\end{equation}
where ${d_0}=(d^{(1)},d^{(2)})^T \in (0,\infty)^2$, ${B_0}=(b_{i,j})_{i,j=1,2}$ is a $2 \times 2$ matrix with non-negative coefficients, and $NB(r,p)$ denotes the negative binomial distribution with parameter $(r,p)$ and mean $r(1-p)/p$. 
%
The parameter $r$ is assumed to be known for each simulation; that is, the parameter of interest is $\theta_0=(d^{(1)},d^{(2)}, b_{1,1},b_{1,2},b_{2,1},b_{2,2})$ and the critical value of the test is $c_{6,\alpha} \approx 4.375$ (see also  Lee \textit{et al.} (2003)). 
 In situations of break, we also assume that the parameter changes from $\theta_0$ (which is characterized here by  $({d_0},{B_0})$) to $\theta_1$ that we will characterize by $({d_1},{B_1})$ for this model. 
\end{itemize}

\begin{figure}[h!]
\begin{center}
\includegraphics[height=5.5cm, width=17cm]{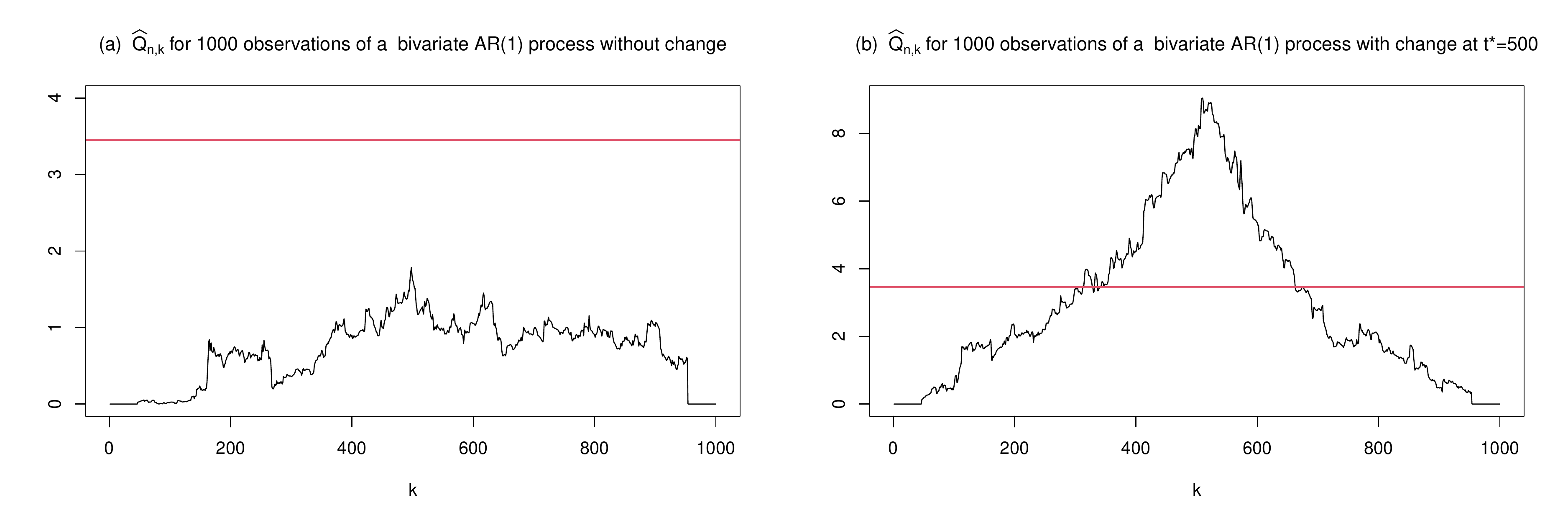} 
\end{center}
\vspace{-.8cm}
 \caption{\small \it Typical realizations of the statistics $\widehat{Q}_{n,k}$  for two trajectories generated from bivariate AR$(1)$ processes.  
(a) is a realization for 1000 observations with a constant parameter $\theta_0=(0.5,-0.2,0.35,0.1)$.
(b) is a realization for 1000 observations in a scenario where the parameter changes from $\theta_0=(0.5,-0.2,0.35,0.1)$  to $\theta_1=(0.5,-0.2,0.1,0.1)$ at $t^*=500$. 
The horizontal line represents the limit of the critical region of the test.}
\label{Graphe_Sim}
\end{figure}

\medskip
\noindent
Figure \ref{Graphe_Sim} is an illustration of a typical realizations of the statistics $\widehat{Q}_{n,k}$ for two trajectories of length $1000$ generated from bivariate AR$(1)$ processes:  a trajectory without change  and a trajectory with a change at time $t^*=500$.
One can see that, for the trajectory without change, the statistics $\widehat{Q}_{n,k}$ are well below the critical value (see Figure \ref{Graphe_Sim}(a)). 
 For the scenario with change, the maximum (which represents the value of $\widehat{Q}_{n}$) of the statistics  $\widehat{Q}_{n,k}$ is higher than the limit of the critical region and that it  is obtained at a point very close to the instant of break (see Figure \ref{Graphe_Sim}(b)).
This empirically comforts the common use of the classical estimator $\widehat{t}_n = \underset{k \in \mathcal T_n} {\text{argmax}} \left(\widehat Q_{n,k}\right)$ to determine the break-point.

\medskip
\noindent 
To evaluate the empirical level and power, 
we consider  trajectories  generated from the two models (\ref{AC_SIM}) and (\ref{MOD_SIM}) in the following situations: 
 (i)  scenarios  with a constant parameter $\theta_0$ and (ii)  scenarios with a parameter change ($\theta_0 \rightarrow \theta_1$)  at time $t^*=n/2$. 
The replication number in each simulation is 200. 
For different scenarios, Table \ref{Table_Sim} indicates the proportion of the number of rejections of the null hypothesis computed under H$_0$ (for the levels) and H$_1$ (for the powers). 
As can be seen from this table, the empirical levels are close to the nominal level for each of the two models. 
%
%
One can see that, the statistic is quite sensitive for detecting the change for both the cases considered under the alternative: the scenario with dependent components and independent components (i.e, the scenario where the matrix $A_1$ or $B_1$ is diagonal) after the breakpoint. 
For both the classes of models, the results of the test are quite accurate; the empirical level approaching the nominal one when $n$ increases and the empirical  power increases with $n$ and is close to 1 when $n=1000$.
This is consistent with the asymptotic results of Theorem \ref{th1_test} and \ref{th2_test}.   

          
 \begin{table}[h!]
\centering
\scriptsize
\caption{\small \it Empirical levels and powers at the nominal level $0.05$ for the change-point test in the models (\ref{AC_SIM}) and (\ref{MOD_SIM}).}
\label{Table_Sim}
\vspace{.3cm}
\begin{tabular}{lllcc}
\hline
\rule[0cm]{0cm}{.5cm}
&&Scenarios&$n=500$&$n=1000$\\
\hline
\rule[0cm]{0cm}{.4cm}
\multirow{8}*{\footnotesize Model (\ref{AC_SIM})}$~$
&Levels: &&&\\
                
                  && $ {A_0}=\begin{pmatrix} 0.6& 0.3\\ 0.0 & 0.4 \end{pmatrix}$                  &0.040& 0.055\\
                  
  &&& &\\
                        
                  &&${A_0}=\begin{pmatrix} 0.5&\!\!\!\!\! -0.2 \\ 0.35 &  0.1 \end{pmatrix}$     &0.060&0.045\\
                
&&&&\\
\rule[0cm]{0cm}{.2cm}
&Powers:  && &\\

                  &&$ {A_0}=\begin{pmatrix} 0.6& 0.3\\ 0.4 & 0.4 \end{pmatrix}~;~ 
                    {A_1}=\begin{pmatrix} 0.6& 0.0\\ 0.0& 0.4 \end{pmatrix}$      &0.765&0.965\\  
 
  &&&&\\

                  &&${A_0}=\begin{pmatrix} 0.5&\!\!\!\!\! -0.2 \\ 0.35 & 0.1 \end{pmatrix}~;~ 
                     {A_1}=\begin{pmatrix} 0.5&\!\!\!\!\!  -0.2 \\ 0.1 & 0.1 \end{pmatrix}$      &0.680&0.940\\

     & &&&\\    
  \hdashline[3pt/3pt] 
   & &&&\\

\multirow{8}*{\footnotesize Model (\ref{MOD_SIM})}$~$   
&Levels: &&&\\
                
                  && ${d_0}=(1,0.3)^{T} \text{ and } {B_0}=\begin{pmatrix} 0.5&0.2\\ 0.1&0.4 \end{pmatrix}$    &0.065&0.055\\
                  
 &&& &\\  
                 
                  &&${d_0}=(0.5,0.5)^{T} \text{ and } {B_0}=\begin{pmatrix} 0.25&0.5 \\ 0.1&0.35 \end{pmatrix}$    &0.065&0.050\\
                
&&&&\\
\rule[0cm]{0cm}{.2cm}
&Powers:  && &\\

                  &&${d_0}=(1,0.3)^{T} \text{ and } {B_0}=\begin{pmatrix} 0.5&0.2\\ 0.1&0.4 \end{pmatrix}~;~ 
                     {d_1}=d_0 \text{ and } {B_1}=\begin{pmatrix} 0.5&0.0\\ 0.0&0.4\end{pmatrix}$     &0.840&0.985\\ 
                     
 &&&&\\

                  &&${d_0}=(0.5,0.5)^{T} \text{ and } {B_0}=\begin{pmatrix} 0.25&0.5 \\ 0.1&0.35 \end{pmatrix}~;~ 
                      {d_1}=(0.5,1)^{T} \text{ and } {B_1}=B_0$     &0.975&0.995\\   
 \Xhline{.75pt}
\end{tabular}
\end{table}



 \section{Proofs of the main results}  
  Let $ (\psi_n)_{n \in \N}  $ and $ (r_n)_{n \in \N}  $ be sequences of random variables or vectors. Throughout this section, we use the notation
  $ \psi_n = o_P(r_n)  $ to mean:  for all $  \varepsilon  > 0, ~ \prob( \|\psi_n \| \geq \varepsilon \|r_n \| ) \limiten 0$.
  Write $ \psi_n = O_P(r_n)$ to mean:  for all  $  \varepsilon > 0 $,  there exists $C>0$  such that  $\prob( \|\psi_n \| \geq C \|r_n \| )\leq \varepsilon $ 
   for $n$ large enough. 
   In the sequel, $C$ denotes a positive constant  whose the value may differ from one inequality to
another.

\subsection{Proof of the results of Section \ref{Sec_detection}}

\subsubsection{Proof of Theorem \ref{th1_test}}
 Define the statistic
 \begin{align*}
  Q_n=\max_{k \in \mathcal T_n} \big(Q_{n,k}\big)
 ~~~\text{with}~~~
  Q_{n,k}=\frac{\left(k(n-k)\right)^{2}}{n^{3}}
 \big(\widehat{\theta}(T_{1,k})-\widehat{\theta}(T_{k+1,n})\big)^T
 \Omega \big(\widehat{\theta}(T_{1,k})-\widehat{\theta}(T_{k+1,n})\big),
 \end{align*}
where $\Omega$ is the covariance matrix defined in the assumption (\textbf{A4}). 
 For any segment $T \subset \{1,\cdots,n\}$ and $\theta \in \Theta$, we also define the function
 \begin{equation*}
C(T,\theta) = \sum_{t \in T} \varphi_t(\theta),~ 
\text{where }  (\varphi_t(\cdot))_{t \in \Z} \text{ is  given in (\textbf{A3})}.
  \end{equation*} 
Let $1\leq k \leq k^{\prime} \leq n$, $\bar{\theta} \in \Theta$ and $i \in \{1,2,\cdots,d\}$.
By the mean value theorem applied to the function $\theta \mapsto \frac{\partial}{\partial \theta_i}C({T_{k,k^{\prime}}},\theta)$, there exists $\theta_{n,i}$ between $\bar{\theta}$ and $\theta^*$ such that
\[ 
\frac{\partial}{\partial \theta_i} C({T_{k,k^{\prime}}},\bar{\theta})=\frac{\partial}{\partial \theta_i} C({T_{k,k^{\prime}}},\theta^*) +\frac{\partial^{2}}{\partial \theta\partial \theta_i} C({T_{k,k^{\prime}}},\theta_{n,i})(\bar{\theta}-\theta^*),
\]
which implies 
\begin{eqnarray}\label{test_Eq_Taylor}
(k^{\prime}-k+1)F_n({T_{k,k^{\prime}}},\bar{\theta}) (\theta^*-\bar{\theta})=\frac{\partial}{\partial \theta} C({T_{k,k^{\prime}}},\theta^*)-\frac{\partial}{\partial \theta} C({T_{k,k^{\prime}}},\bar{\theta})
\end{eqnarray}
with
\begin{equation}\label{test_Def_F_n}
F_n(T_{k,k^{\prime}},\bar{\theta})= \frac{1}{(k^{\prime}-k+1)}\frac{\partial^{2}}{\partial \theta\partial \theta_i} C({T_{k,k^{\prime}}},\theta_{n,i})_{1\leq i \leq d}.
\end{equation}  
\medskip

\medskip
\noindent 
The following lemma will be useful in the sequel.
\begin{lem}\label{test_Lem1}
Assume that the conditions  of Theorem \ref{th1_test} hold. 
\begin{enumerate}
\rm
    \item [(i)] $\underset{k \in \mathcal T_n}{\max}\big|\widehat Q_{n,k}-Q_{n,k}\big|=o_P(1)$.
    
    
\rm
\item [(ii)]  \it
If $(j_n)_{n\geq1}$ and $(k_n)_{n\geq1}$ are two integer valued sequences such that $j_n \leq k_n$, $k_n \limiten \infty$ and $k_n - j_n \limiten \infty $, then 
$
F_n(T_{j_n,k_n},\widehat{\theta}(T_{j_n,k_n})) \limitepsn F, 
$
 where  $F$  is the matrix defined in   (\textbf{A4}).
\end{enumerate}
\end{lem}

\medskip
\noindent
{\bf Proof.}
\begin{enumerate}
 \item [(i)] Let $k \in \mathcal T_n$. 
As $n \rightarrow \infty$,   
from the asymptotic normality of the MCE  and the consistency of $\widehat{\Omega}(u_n)$, we obtain:
\begin{equation}\label{cov_OP_H0}
\big\|\sqrt{k} \big(\widehat{\theta}(T_{1,k}) -\theta^*\big)\big\|=O_P(1), ~~
\big\|\sqrt{n-k} \big(\widehat{\theta}(T_{k+1,n}) -\theta^*\big)\big\|=O_P(1) ~~\text{and}~~\big\| \widehat{\Omega}(u_n)-\Omega\big\|=o(1).
\end{equation}
Then, it holds that
\begin{align*}
\big|\widehat Q_{n,k}-Q_{n,k}\big|
&=
\frac{\left(k(n-k)\right)^{2}}{n^{3}}
 \Big|
\big(\widehat{\theta}(T_{1,k})-\widehat{\theta}(T_{k+1,n})\big)^T 
 \big( \widehat{\Omega}(u_n)-\Omega\big)
\big(\widehat{\theta}(T_{1,k})-\widehat{\theta}(T_{k+1,n})\big)
\Big|\\
&\leq
C\frac{\left(k(n-k)\right)^2}{n^{3}}
\big\| \widehat{\Omega}(u_n)-\Omega\big\|
\big\|\widehat{\theta}(T_{1,k})-\widehat{\theta}(T_{k+1,n}) \big\|^2\\
&\leq
 C\big\| \widehat{\Omega}(u_n)-\Omega\big\| 
 \bigg[ 
\frac{k(n-k)^2}{n^{3}} \big\|\sqrt{k} \big(\widehat{\theta}(T_{1,k}) -\theta^*\big)\big\|^2 
+
\frac{k^2(n-k)}{n^{3}} \big\|\sqrt{n-k} \big(\widehat{\theta}(T_{k+1,n}) -\theta^*\big)\big\|^2 
 \bigg]\\
&\leq
 o(1)O_P(1);
\end{align*}
which allows to conclude.
 

\item [(ii)] Applying $(\ref{test_Def_F_n})$ with $\bar{\theta} = \widehat{\theta}(T_{j_n,k_n})$, we obtain
\begin{align*}
F_n(T_{j_n,k_n},\widehat{\theta}(T_{j_n,k_n}))
&=  \Big(\frac{1}{k_n-j_n+1}\frac{\partial^{2}}{\partial \theta\partial \theta_i} C(T_{j_n,k_n},\theta_{n,i})\Big)_{1\leq i \leq d}\\
&=\frac{1}{k_n-j_n+1}\Big(\sum_{t \in T_{j_n,k_n}}\frac{\partial^{2} \varphi_t(\theta_{n,i})}{\partial \theta\partial \theta_i} \Big)_{1\leq i \leq d},
\end{align*}
where $\theta_{n,i}$ belongs between $\widehat{\theta}(T_{j_n,k_n})$ and $\theta^*$.  
Since $\widehat{\theta}(T_{j_n,k_n})  \limitepsn \theta^*$, $~ \theta_{n,i} ~ \limitepsn \theta^*$ (for any  $i=1,\cdots,d$)
 and that $F=\E\big[\frac{\partial^{2} \varphi_0(\theta^*)}{\partial \theta\partial\theta^T} \big]$ exists (see the assumption (\textbf{A4})), by the uniform strong law of large numbers, for any $i =1,\cdots,d$, we get 
\begin{align*}
 & \Big \| \frac{1}{k_n-j_n+1}\sum_{t \in T_{j_n,k_n}}\frac{\partial^{2}  \varphi_t(\theta_{n,i})}{\partial \theta\partial \theta_i} 
-  \E\Big[\frac{\partial^{2} \varphi_0(\theta^*)}{\partial \theta\partial\theta_i}\Big] \Big \| \\
&\leq \Big \| \frac{1}{k_n-j_n+1}\sum_{t \in T_{j_n,k_n}}\frac{\partial^{2}  \varphi_t(\theta_{n,i})}{\partial \theta\partial \theta_i} 
-  \E\Big[\frac{\partial^{2} \varphi_0(\theta_{n,i})}{\partial \theta\partial\theta_i}\Big] \Big \| 
+  \Big \| \E\Big[\frac{\partial^{2} \varphi_0(\theta_{n,i})}{\partial \theta\partial\theta_i}\Big] 
-  \E\Big[\frac{\partial^{2} \varphi_0(\theta^*)}{\partial \theta\partial\theta_i}\Big] \Big \|  \\
&\leq \Big \| \frac{1}{k_n-j_n+1}\sum_{t \in T_{j_n,k_n}}\frac{\partial^{2}  \varphi_t(\theta)}{\partial \theta\partial \theta_i} 
-  \E\Big[\frac{\partial^{2} \varphi_0(\theta)}{\partial \theta\partial\theta_i}\Big] \Big \|_\Theta 
+  o(1) = o(1) + o(1)= o(1).
\end{align*}
%
This completes the proof of the lemma.
\begin{flushright}
$\blacksquare$ 
\end{flushright}
\end{enumerate} 

\noindent 
Now, we use (\ref{cov_OP_H0}) and the part (ii) of Lemma \ref{test_Lem1} to show that
\begin{equation}\label{Cond_proof_th1}
Q_n \limiteloin \sup_{0\leq \tau\leq 1}\left\|W_d(\tau)\right\|^{2}.
\end{equation}
\noindent Let $k \in \mathcal T_n$. Applying (\ref{test_Eq_Taylor}) with $\bar{\theta}=\widehat{\theta}(T_{1,k})$ and $T_{k,k^\prime}=T_{1,k}$, we get
\begin{eqnarray}\label{test_Eq_T_1,k}
F_n(T_{1,k},\widehat{\theta}(T_{1,k})) \cdot (\theta^*-\widehat{\theta}(T_{1,k}))=\frac{1}{k}\Big(\frac{\partial}{\partial \theta} C(T_{1,k},\theta_1^*)-\frac{\partial}{\partial \theta} C(T_{1,k},\widehat{\theta}(T_{1,k}))\Big).
\end{eqnarray}
With $\bar{\theta}=\widehat{\theta}(T_{k+1,n})$ and $T_{k,k^\prime}=T_{k+1,n}$, (\ref{test_Eq_Taylor}) becomes
\begin{eqnarray}\label{test_Eq_T_k+1,n}
F_n(T_{k+1,n},\widehat{\theta}(T_{k+1,n})) \cdot (\theta^*-\widehat{\theta}(T_{k+1,n}))=\frac{1}{n-k}\Big(\frac{\partial}{\partial \theta} C(T_{k+1,n},\theta_1^*)-\frac{\partial}{\partial \theta} C(T_{k+1,n},\widehat{\theta}(T_{k+1,n}))\Big).
\end{eqnarray}
Moreover, as $n\rightarrow +\infty$, Lemma \ref{test_Lem1}(ii) implies
\begin{align*}
\big\|F_n(T_{1,k},\widehat{\theta}(T_{1,k}))-F\big\|=o(1) 
~\text{ and }~
\big\|F_n(T_{k+1,n},\widehat{\theta}(T_{k+1,n}))-F\big\|=o(1).
\end{align*}
Then, according to (\ref{cov_OP_H0}), for  $n$ large enough, (\ref{test_Eq_T_1,k}) gives 
\begin{align*}
&\sqrt{k}F\left(\theta^*-\widehat{\theta}(T_{1,k})\right)\\
&\hspace{.5cm}=
\frac{1}{\sqrt{k}}\Big(\frac{\partial}{\partial \theta} C(T_{1,k},\theta^*)-\frac{\partial}{\partial \theta} C(T_{1,k},\widehat{\theta}(T_{1,k}))\Big) 
-\sqrt{k}\left(\left(F_n(T_{1,k},\widehat{\theta}(T_{1,k}))-J\right)\left(\widehat{\theta}(T_{1,k})-\theta_0\right)\right)\\
&\hspace{.5cm}= \frac{1}{\sqrt{k}}\Big(\frac{\partial}{\partial \theta} C(T_{1,k},\theta^*)-\frac{\partial}{\partial \theta} C(T_{1,k},\widehat{\theta}(T_{1,k}))\Big) +o_P(1)\\
&\hspace{.5cm}= \frac{1}{\sqrt{k}}\Big(\frac{\partial}{\partial \theta} C(T_{1,k},\theta^*)- \frac{\partial}{\partial \theta} \widehat C(T_{1,k},\widehat{\theta}(T_{1,k}))\Big)+o_P(1)
+\frac{1}{\sqrt{k}}\Big(\frac{\partial}{\partial \theta} \widehat C(T_{1,k},\widehat{\theta}(T_{1,k}))-\frac{\partial}{\partial \theta} C(T_{1,k},\widehat{\theta}(T_{1,k}))\Big) \\
&\hspace{.5cm}= \frac{1}{\sqrt{k}}\Big(\frac{\partial}{\partial \theta} C(T_{1,k},\theta^*)-\frac{\partial}{\partial \theta} \widehat C(T_{1,k},\widehat{\theta}(T_{1,k}))\Big) +o_P(1)~ ~(\text{from the condition  (\ref{eq_assump_app_dif}) in (\textbf{A3})}) .
\end{align*}
This is equivalent to
\begin{equation}\label{test_Eq_a}
F\left(\theta^*-\widehat{\theta}(T_{1,k})\right)=\frac{1}{k}\Big(\frac{\partial}{\partial \theta} C(T_{1,k},\theta^*)-\frac{\partial}{\partial \theta} \widehat C(T_{1,k},\widehat{\theta}(T_{1,k}))\Big) +o_P\Big(\frac{1}{\sqrt{k}}\Big).
\end{equation}
 For $n$ large enough, $\widehat{\theta}(T_{1,k})$ is an interior point of $\Theta$ and we have $ \frac{\partial}{\partial \theta} \widehat C(T_{1,k},\widehat{\theta}(T_{1,k}))=0$.
Hence, for $n$ large enough, we get from (\ref{test_Eq_a})
\begin{equation}\label{test_Eq_a_bis}
F\left(\theta^*-\widehat{\theta}(T_{1,k})\right)
 =\frac{1}{k}\frac{\partial}{\partial \theta} C(T_{1,k},\theta^*)+o_P\Big(\frac{1}{\sqrt{k}}\Big).
\end{equation}
Similarly, we can use (\ref{test_Eq_T_k+1,n}) to obtain
\begin{equation}\label{test_Eq_a_bisbis}
F\left(\theta^*-\widehat{\theta}(T_{k+1,n})\right)=\frac{1}{n-k}\frac{\partial}{\partial \theta} C(T_{k+1,n},\theta^*)+o_P\Big(\frac{1}{\sqrt{n-k}}\Big).
\end{equation}
The subtraction of  (\ref{test_Eq_a_bis}) and (\ref{test_Eq_a_bisbis})  gives 
\begin{align*}
-F\left(\widehat{\theta}(T_{1,k})-\widehat{\theta}(T_{k+1,n})\right)&=\frac{1}{k}\frac{\partial}{\partial \theta} C(T_{1,k},\theta^*)-\frac{1}{n-k}\frac{\partial}{\partial \theta} C(T_{k+1,n},\theta^*)+o_P\Big(\frac{1}{\sqrt{k}}+\frac{1}{\sqrt{n-k}}\Big)\\
&=\frac{1}{k}\frac{\partial}{\partial \theta} C(T_{1,k},\theta^*)-\frac{1}{n-k}\Big(\frac{\partial}{\partial \theta} C(T_{1,n},\theta^*)-\frac{\partial}{\partial \theta} C(T_{1,k},\theta^*)\Big)
 +o_P\Big(\frac{1}{\sqrt{k}}+\frac{1}{\sqrt{n-k}}\Big)\\
&=\frac{n}{k(n-k)}\Big(\frac{\partial}{\partial \theta} C(T_{1,k},\theta^*)-\frac{k}{n}\cdot\frac{\partial}{\partial \theta} C(T_{1,n},\theta^*)\Big)+o_P\Big(\frac{1}{\sqrt{k}}+\frac{1}{\sqrt{n-k}}\Big).
\end{align*}
Since the matrix $G$ is positive definite (see (\textbf{A3})), the above equality is equivalent to 
\begin{align}\label{test_Eq_b}
&-\frac{k(n-k)}{n^{3/2}}G^{-1/2}F\big(\widehat{\theta}(T_{1,k})-\widehat{\theta}(T_{k+1,n})\big)\nonumber\\
&\hspace{3cm}=
\frac{G^{-1/2}}{\sqrt{n}}\Big(\frac{\partial}{\partial \theta} C(T_{1,k},\theta^*)-\frac{k}{n}\cdot\frac{\partial}{\partial \theta} C(T_{1,n},\theta^*)\Big) +o_P\Big(\frac{\sqrt{k(n-k)}}{n} + \frac{ \sqrt{n-k}}{\sqrt{n}}\Big) \nonumber\\
&\hspace{3cm}=\frac{G^{-1/2}}{\sqrt{n}}\Big(\frac{\partial}{\partial \theta} C(T_{1,k},\theta^*)-\frac{k}{n}\cdot\frac{\partial}{\partial \theta} C(T_{1,n},\theta^*)\Big) +o_P(1)
\end{align}
and that $Q_{n,k}$ can be rewritten as 
\begin{equation}\label{test_Def2_Q_{n,k}}
Q_{n,k}=\big\|\frac{k(n-k)}{n^{3/2}}G^{-1/2}F \big(\widehat{\theta}(T_{1,k})-\widehat{\theta}(T_{k+1,n})\big)\big\|^2~ \text{ for all } k \in \mathcal T_n.
\end{equation}
Moreover, applying the central limit theorem for the martingale difference sequence $\left(\frac{\partial}{\partial \theta} \varphi_t(\theta^*),\mathcal{F}_{t}\right)_{t \in \mathbb{Z}}$ (see \nocite{Billingsley1968} Billingsley (1968)), we have
\begin{align*}
\frac{1}{\sqrt{n}}\Big(\frac{\partial}{\partial \theta}C(T_{1,[n \tau]},\theta^*)-\frac{[n \tau]}{n}\frac{\partial}{\partial \theta}C(T_{1,n },\theta^*)\Big)
&=
\frac{1}{\sqrt{n}}\Big(\sum_{t=1}^{[n \tau]}\frac{\partial}{\partial \theta}\varphi_t(\theta^*)-\frac{[n \tau]}{n}\sum_{t=1}^{n }\frac{\partial}{\partial \theta}\varphi_t(\theta^*)\Big)\\
&~~
\limiteloin B_{G}(\tau)-\tau B_{G}(1),
\end{align*}
where $[x]$  denotes the integer part  of $x$ and  $B_{G}$ is a Gaussian process with covariance matrix $\min(s,t)G$. \\
Then,
\begin{align*}
\frac{G^{-1/2}}{\sqrt{n}}\Big(\frac{\partial}{\partial \theta}C(T_{1,[n \tau]},\theta^*)-\frac{[n \tau]}{n}\frac{\partial}{\partial \theta}C(T_{1,n },\theta^*)\Big)
\limiteloin 
B_{d}(\tau)-\tau B_{d}(1)=W_d(\tau)  ~\text{ in }  D([0,1]),
\end{align*}
 where $B_d$ is a $d$-dimensional standard motion and $W_d$ is a $d$-dimensional Brownian bridge.\\
Therefore, using (\ref{test_Eq_b}) and (\ref{test_Def2_Q_{n,k}}), we obtain
\begin{align*}
Q_{n,[n\tau]}&=
\Big\|\frac{[n\tau](n-[n\tau])}{n^{3/2}}G^{-1/2}F\big(\widehat{\theta}(T_{1,[n\tau]})-\widehat{\theta}(T_{[n\tau]+1,n})\big)\Big\|^2
\limiteloin 
\sup_{0\leq \tau \leq  1} \left\|W_d(\tau)\right\|^2~\text{ in }~ D([0,1]).
\end{align*}
For $n$ large enough, we deduce
\[
Q_n=\max_{v_n \leq k  \leq n-v_n} \big(Q_{n,k}\big)=
\sup_{\frac{v_n}{n}\leq \tau \leq 1-\frac{v_n}{n}} Q_{n,[n\tau]} \limiteloin 
\sup_{0\leq \tau \leq  1} \left\|W_d(\tau)\right\|^2 ~\text{ in }~ D([0,1]);\]
which shows that (\ref{Cond_proof_th1}) holds. 
Hence, we can conclude the proof of the theorem from Lemma \ref{test_Lem1}(i).
\begin{flushright}
$\blacksquare$ 
\end{flushright}

\subsubsection{Proof of Theorem \ref{th2_test}}
Under the alternative, we can write 
\begin{equation*}
Y_{t}=\left\{
\begin{array}{ll}
Y^{(1)}_{t}~~\textrm{for}~~t \leq t^*,\\
\rule[0cm]{0cm}{.7cm}
Y^{(2)}_{t}~~\textrm{for}~~t >t^*,
\end{array}
\right.
\end{equation*} 
where $t^*=[\tau^* n]$ (with $0<\tau^*<1$) and $\{Y^{(j)}_{t}, t \in \mathbb{Z}\}$ ($j=1,2$) is a stationary and ergodic process depending on the parameter $\theta^{*}_j$ (with $\theta^{*}_1 \neq \theta^{*}_2$) satisfying the assumptions (\textbf{A1})-(\textbf{A4}).\\
Remark that $ \widehat{Q}_n = \underset{k \in \mathcal T_n}{\max}\big(\widehat{Q}_{n,k}\big) \geq \widehat{Q}_{n,t^*} $.  Then,
to prove the theorem, we will show that $\widehat{Q}_{n,t^*} \limiteproban +\infty$.

\medskip
\noindent
For any $n\in \N$, we have
\[
 \widehat{Q}_{n,t^{*}}
 =
\frac{\left(t^*(n-t^{*})\right)^2}{n^{3}}\big(\widehat{\theta}(T_{1,t^{*}})-\widehat{\theta}(T_{t^{*}+1,n})\big)^T\widehat{\Omega}(u_n)\big(\widehat{\theta}(T_{1,t^{*}})-\widehat{\theta}(T_{t^{*}+1,n})\big)
\]
with
\[
\widehat{\Omega}(u_n)
=\frac{1}{2} 
\left[
\widehat F(T_{1,u_n})  \widehat G(T_{1,u_n})^{-1}   \widehat F(T_{1,u_n}) +
\widehat F(T_{u_n+1,n})  \widehat G(T_{u_n+1,n})^{-1}  \widehat F(T_{u_n+1,n})
\right].
\]
%
Moreover, the two matrices in the formula of $\widehat{\Omega}(u_n)$ are  positive
semi-definite.
Then, we obtain
\begin{align}
 \widehat{Q}_{n,t^{*}}
&= 
\frac{\left([\tau^* n](n-[\tau^* n])\right)^2}{n^{3}}
 \left(\widehat{\theta}(T_{1,t^{*}})-\widehat{\theta}(T_{t^{*}+1,n})\right)^T \nonumber\\
& \hspace{4cm} \times \Big[
\widehat F(T_{1,u_n})  \widehat G(T_{1,u_n})^{-1}   \widehat F(T_{1,u_n}) +
\widehat F(T_{u_n+1,n})  \widehat G(T_{u_n+1,n})^{-1}  \widehat F(T_{u_n+1,n})
\Big]\nonumber\\
& \hspace{5cm} \times \left(\widehat{\theta}(T_{1,t^{*}})-\widehat{\theta}(T_{t^{*}+1,n})\right)\nonumber\\
&\label{Major_Q_n_t}\geq 
 n
\left(\widehat{\theta}(T_{1,t^{*}})-\widehat{\theta}(T_{t^{*}+1,n})\right)^T
\Big[
\widehat F(T_{1,u_n})  \widehat G(T_{1,u_n})^{-1}   \widehat F(T_{1,u_n}) 
\Big]
\left(\widehat{\theta}(T_{1,t^{*}})-\widehat{\theta}(T_{t^{*}+1,n})\right).
\end{align}
By the consistency and asymptotic normality of the MCE, we have: (i)
$
 \widehat{\theta}(T_{1,t^{*}})-\widehat{\theta}(T_{t^*+1,n}) \limitepsn \theta^{*}_1-\theta^{*}_{2}\neq 0
 $
and (ii) 
$\widehat F(T_{1,u_n})  \widehat G(T_{1,u_n})^{-1}   \widehat F(T_{1,u_n}) $
converges to the covariance
matrix of the stationary model of the first regime which is positive definite. 
Therefore,  (\ref{Major_Q_n_t}) implies $\widehat{Q}_{n,t^*} \limitepsn +\infty$. 
This establishes the theorem.
\begin{flushright}
$\blacksquare$ 
\end{flushright}

 \subsection{Proof of the results of Section \ref{Sec_class_MOD}}
  \subsection{Proof of Proposition \ref{prop1}}
 Let $F_\lambda(y)$  be the cumulative distribution function of $p(y|\eta)$ with marginals $F_{\lambda_{1},1}, \cdots, F_{\lambda_{m},m}$, where $\lambda = (\lambda_{1},\cdots,\lambda_{m})^T = \partial A (\eta)/\partial \eta$. 
From the Sklar's theorem (see \nocite{Sklar1959} Sklar (1959)), one can find a copula $\mathcal{C}$ such that, for all $y=(y_1,\cdots, y_m) \in \R^m$
\[ F_\lambda(y) = \mathcal{C} \big(F_{\lambda_{1},1}(y_1),\cdots,F_{\lambda_{m},m}(y_m) \big) .\]
 For $i=1,\cdots,m$, denote by $F_{\lambda,i}^{-1}(u) := \inf\{ y_i \geq 0, ~  F_{\lambda,i}(y_i) \geq u \} $ for all $u \in [0,1]$.
Let $\{ U_t = (U_{t,1},\cdots,U_{t,m})^T, ~ t \in \Z \}$ be a sequence of independent  random vectors with distribution $\mathcal{C}$. 
 We will prove that, there exists  a $\tau$-weakly
dependent, stationary and ergodic solution $(Y_t, \lambda_t)$ of (\ref{Model_exp_f}) satisfying:
\begin{equation}\label{Y_t_sol}
 Y_t = \big(F_{\lambda_{t,1},1}^{-1}(U_{t,1}),\cdots,F_{\lambda_{t,m},m}^{-1}(U_{t,m}) \big)^T
\end{equation}
with $\lambda_t = \left(\lambda_{t,1},\cdots,\lambda_{t,m} \right)^T= f_\theta(Y_{t-1},\cdots)$.
For a process $(Y_t)_{t \in \Z}$ that fulfills (\ref{Model_exp_f}) and (\ref{Y_t_sol}), we get,
\begin{equation}\label{Y_t_Psi}
Y_t = \big(F_{\lambda_{t,1},1}^{-1}(U_{t,1}),\cdots,F_{\lambda_{t,m},m}^{-1}(U_{t,m}) \big)^T :=\Psi(Y_{t-1},\cdots; U_t),
\end{equation}
where $\Psi$ is a function defined in $(\N_0^m)^\infty \times [0,1]^m$. 
According to \nocite{Doukhan2008} Doukhan and Wintenberger (2008), it suffices to show that: (i)  $\E \| \Psi(\pmb{y};U_t) \| < \infty $ for some $\pmb{y} \in (\N_0^m)^\infty$  and (ii) there exists a sequence of non-negative real numbers $(\alpha_k(\Psi))_{k \geq 1}$ satisfying $\sum_{k \geq 1} \alpha_k(\Psi) < 1$ such that, for all $\pmb{y}, \pmb{y}' \in (\N_0^m)^\infty$,
 $\E \| \Psi(\pmb{y};U_t) - \Psi(\pmb{y}';U_t) \| \leq \sum_{k \geq 1} \alpha_k(\Psi) \|y_k - y'_k \| $. 
  \medskip
  
 \textbf{Proof of (i)}.
 Set $f_\theta(0,\cdots) = \lambda = (\lambda_{1},\cdots,\lambda_{m})^T$. The random vector $\big(F_{\lambda_{1},1}^{-1}(U_{t,1},\cdots,F_{\lambda_{m},m}^{-1}(U_{t,m}) \big)^T$ is $F_\lambda$ distributed.
Thus,
\[ \E \| \Psi(0,\cdots;U_t) \| = \E \big \| \big(F_{\lambda_{1},1}^{-1}(U_{t,1}),\cdots,F_{\lambda_{m},m}^{-1}(U_{t,m}) \big)  \big \| = \| \lambda \| = \| f_\theta(0,\cdots) \| < \infty,  \] 
 where this inequality holds from the assumption \textbf{A}$_0(\Theta)$.
 
  \medskip
  
 \textbf{Proof of (ii)}.
  For all $\pmb{y}, \pmb{y}' \in (\N_0^m)^\infty$, set $\lambda = f_\theta(\pmb{y},\cdots) = (\lambda_{1},\cdots,\lambda_{m})^T$ and $\lambda' = f_\theta(\pmb{y}',\cdots) = (\lambda_{1}',\cdots,\lambda_{m}')^T$. We have,
 \begin{align}
  \nonumber \E \| \Psi(\pmb{y};U_t) - \Psi(\pmb{y}';U_t) \|
   &= 
    \E \big \| \big(F_{\lambda_{1},1}^{-1}(U_{t,1}),\cdots,F_{\lambda_{m},m}^{-1}(U_{t,m}) \big)  -  \big(F_{\lambda_{1}',1}^{-1}(U_{t,1}),\cdots,F_{\lambda_{m}',m}^{-1}(U_{t,m}) \big)  \big\| \\
 \label{Lip_F_lambda} 
 &\leq 
 \sum_{i=1}^m \E \big| F_{\lambda_{i},i}^{-1}(U_{t,i}) - F_{\lambda_{i}',i}^{-1}(U_{t,i})  \big|
  =  \sum_{i=1}^m | \lambda_{i}  - \lambda_{i}'| \\
 \label{Lip_f_alpha} & =\| \lambda - \lambda' \| = \| f_\theta(\pmb{y},\cdots) - f_\theta(\pmb{y}',\cdots) \| \leq \sum_{k=1}^{\infty}\alpha^{(0)}_k \|y_k-y'_k\|, 
 \end{align}   
where the equality in (\ref{Lip_F_lambda}) holds from the  Proposition A.2 of \nocite{Davis2016} Davis and Liu (2016) and the inequality in (\ref{Lip_f_alpha}) holds from the assumption \textbf{A}$_0(\Theta)$. 
Thus, take $\alpha_k(\Psi) = \alpha^{(0)}_k$, which completes the proof of the proposition.  

\begin{flushright}
$\blacksquare$
\end{flushright}

 \subsubsection{Proof of Theorem \ref{th1}}
 
%
%
 To simplify, we will use the following notations in the sequel: 
\begin{align*}
\ell_{t,i}(\theta) &:= Y_{t,i}\log \lambda_{t,i}(\theta)- \lambda_{t,i}(\theta) = Y_{t,i}\log f_\theta^{t,i}- f_\theta^{t,i},\\
\widehat \ell_{t,i}(\theta) &:= Y_{t,i}\log \widehat \lambda_{t,i}(\theta)- \widehat \lambda_{t,i}(\theta) = Y_{t,i}\log \widehat f_\theta^{t,i}- \widehat f_\theta^{t,i}
,
\end{align*}
where $ f_\theta^{t,i}$ and $\widehat f_\theta^{t,i}$ (for $i=1,\cdots,m$) represent the $i$th component of
 $f^t_\theta \equiv f_\theta(Y_{t-1}, Y_{t-2}, \cdots )$ and $\widehat  f^t_\theta \equiv f_\theta(Y_{t-1}, Y_{t-2}, \cdots,Y_1 )$, respectively.
 
  \medskip

 (i) Firstly, we will show that 
\begin{equation}\label{appro_Ln}
\frac{1}{n}\big\|\widehat L_{n}(\theta) - L_{n}(\theta)  \big\|_\Theta \limitepsn 0 .
 \end{equation}
%
Remark that
\begin{align}\label{Ln_hat.Ln}  
         \frac{1}{n}\big\|\widehat L_{n}(\theta) - L_{n}(\theta)  \big\|_\Theta
         &\leq 
 \frac{1}{n}\sum_{t=1}^{n}\|\widehat{\ell}_{t}(\theta)-\ell_{t}(\theta) \|_\Theta
 \leq 
 \frac{1}{n}\sum_{i=1}^{m}\sum_{t=1}^{n}\|\widehat{\ell}_{t,i}(\theta)-\ell_{t,i}(\theta) \|_\Theta .  
 \end{align} 
Using $\textbf{A}_0(\Theta)$ with the condition (\ref{cond_th1}) and the existence of the moment of order 2 (i.e., (\ref{moment}) with $\epsilon \geq 1$), one can proceed as in the proof of Theorem 3.1 in \cite{Doukhan2015} to prove that 
 \[
 \frac{1}{n}\sum_{t=1}^{n}\|\widehat{\ell}_{t,i}(\theta)-\ell_{t,i}(\theta) \|_\Theta \limitepsn 0 \text{ for all } i=1,\cdots,m.
 \]
Therefore, (\ref{appro_Ln}) is obtained by using (\ref{Ln_hat.Ln}). 

\medskip

\medskip

(ii) Let us establish that: for all $t \in \Z$,
 \begin{equation}\label{ell_t_finite}
 \E\big[ \|\ell_t(\theta) \|_{\Theta} \big] < \infty .
\end{equation}  
We have 
\[ \E\big[ \|\ell_t(\theta) \|_{\Theta} \big] \leq \sum_{i=1}^m \E\big[ \sup_{\theta\in \Theta} |\ell_{t,i}(\theta) | \big], \]
%
%
From $\textbf{A}_0(\Theta)$, ($\mathcal{MOD}.\textbf{A0}$), (\ref{moment}) (with $\epsilon \leq 1$) and by going along similar lines as in the proof of Theorem 3.1 in \cite{Doukhan2015}, we get: 
$\E\big[ \sup_{\theta\in \Theta} |\ell_{t,i}(\theta) | \big] < \infty$ for all $i=1,\cdots,m$.
Thus, (\ref{ell_t_finite}) holds.

\medskip
\noindent 
Since $\{Y_{t},~t\in \Z\}$ is stationary and ergodic, 
 the process $\{\ell_{t}(\theta),~t\in \Z\}$ is also a stationary and ergodic sequence. Then, by the uniform strong law of large numbers applied to $\{\ell_{t}(\theta),~t\in \Z\}$, it holds that
\begin{equation*} 
\Big\| \frac{1}{n}L_{n}(\theta)-\E [\ell_{0}(\theta)]\Big\|_{\Theta}= \Big\|  \frac{1}{n} \sum_{t=1}^{n} \ell_{t}(\theta)-\E [\ell_{0}(\theta)]\Big\|_{\Theta}  \limitepsn 0.
 \end{equation*}
 Thus,  from (\ref{appro_Ln}), we obtain 
\begin{equation*} 
\Big\| \frac{1}{n}\widehat{L}_{n}(\theta)-\E [\ell_{0}(\theta)]\Big\|_{\Theta} \leq
\frac{1}{n}\big\| \widehat{L}_{n}( \theta)-L_{n}(\theta)\big\|_{\Theta} +\Big\| \frac{1}{n}L_{n}(\theta)-\E [\ell_{0}(\theta)]\Big\|_{\Theta} 
  \limitepsn 0.
 \end{equation*}

\medskip

\medskip
 (iii) To complete the proof of the theorem, it suffices to show that, the function $\theta \mapsto L(\theta) = \E [\ell_0 (\theta)]$ has a unique maximum at $\theta^*$.
Let $\theta \in \Theta$, such that $\theta \neq \theta^*$. We have
\begin{align*}
  L(\theta^*) - L(\theta) &= \sum_{i=1}^m \big( \E \ell_{0,i}(\theta^*) - \E \ell_{0,i}(\theta) \big) = \sum_{i=1}^m \big( \E[Y_{0,i}\log f_{\theta^*}^{0,i} - f_{\theta^*}^{0,i}] -  \E[Y_{0,i}\log f_{\theta}^{0,i} - f_{\theta}^{0,i}] \big) \\
&= \sum_{i=1}^m \big( \E[ f_{\theta^*}^{0,i} \log f_{\theta^*}^{0,i} - f_{\theta^*}^{0,i}] -  \E[ f_{\theta^*}^{0,i} \log f_{\theta}^{0,i} - f_{\theta}^{0,i}] \big)\\
&= \sum_{i=1}^m \Big( \E \big[ f_{\theta^*}^{0,i} \big(\log  f_{\theta^*}^{0,i} - \log f_{\theta}^{0,i} \big) \big] -  \E(f_{\theta^*}^{0,i} - f_{\theta}^{0,i}) \Big).
  \end{align*}
According to the identifiability assumption $\textbf{A}_0(\Theta)$ and since $\theta \neq \theta^*$, there exists $i_0$ such that $f_{\theta}^{0,i_0} \neq f_{\theta^*}^{0,i_0}$.  
By going as in the proof of Theorem 3.1 in \cite{Doukhan2015}, we get $\E \big[ f_{\theta^*}^{0,i_0} \big(\log  f_{\theta^*}^{0,i_0} - \log f_{\theta}^{0,i_0} \big) \big] -  \E[f_{\theta^*}^{0,i_0} - f_{\theta}^{0,i_0}] > 0$  and
$\E \big[ f_{\theta^*}^{0,i} \big(\log  f_{\theta^*}^{0,i} - \log f_{\theta}^{0,i} \big) \big] -  \E[f_{\theta^*}^{0,i} - f_{\theta}^{0,i}] \geq 0$ for $i=1,\cdots,m, ~ i \neq i_0$. This establishes (iii); which consequently yields the theorem.
\begin{flushright}
$\blacksquare$
\end{flushright}
 
 %
%
 \subsubsection{Proof of Theorem \ref{th2}}
  Applying the mean value theorem to the function $\theta \mapsto \frac{\partial}{\partial \theta_i} L_n (\theta)$  for all $i \in \{1,\cdots,d\}$, there exists 
	$\bar  \theta_{n,i}$ between $\widehat{\theta}_n$ and  $\theta^*$ such that
 \begin{equation*}\label{mvt_dl}
 \frac{\partial}{\partial \theta_i}L_n ( \widehat{\theta}_n)= \frac{\partial }{\partial \theta_i}L_n (\theta^*) + \frac{\partial^2}{\partial \theta \partial \theta_i}L_n (\bar \theta_{n,i}) (\widehat{\theta}_n-\theta^*),
\end{equation*}
which is equivalent to 
\begin{equation}\label{mvt.1_dl}
\sqrt{n}J(\widehat{\theta}_n) (\widehat{\theta}_n-\theta^*) = 
\frac{1}{\sqrt{n}} \Big(\frac{\partial}{\partial \theta} L_n (\theta^*) -\frac{\partial}{\partial \theta}\widehat L_n ( \widehat{\theta}_n)\Big) 
+
\frac{1}{\sqrt{n}} \Big( \frac{\partial}{\partial \theta}\widehat L_n (\widehat{\theta}_n)-\frac{\partial}{\partial \theta}L_n (\widehat{\theta}_n)\Big)
\end{equation}
with
\begin{equation}\label{def_Jn}
J(\widehat{\theta}_n) = \big(-\frac{1}{n}\frac{\partial^2}{\partial \theta \partial \theta_i} L_n (\bar \theta_{n,i})\big)_{1 \leq i \leq d}. 
\end{equation}

 \medskip
 \noindent 
The following lemma is needed.
 \begin{lem}\label{lem2} 
 Assume that the conditions of Theorem \ref{th2} hold. Then,
 \begin{enumerate}
 \item[(i)] $ \E\Big[ \frac{1}{\sqrt{n}} \Big\|\frac{\partial}{\partial \theta}\widehat L_n (\theta)- \frac{\partial}{\partial \theta} L_n (\theta)\Big\|_\Theta\Big]  \limiten 0. $
 \item[(ii)] $ \frac{1}{n} \Big\|\frac{\partial^2}{\partial \theta \partial \theta^T}\widehat L_n (\theta)- \frac{\partial^2}{\partial \theta \partial \theta^T} L_n (\theta)\Big\|_\Theta = o(1).$
 \end{enumerate}
  \end{lem}
  
  \medskip
  \noindent
 \emph{\bf Proof.}\\
(i) We have
\begin{align}\label{dif_Ln_hat.Ln}  
         \frac{1}{\sqrt{n}}  \Big\| \frac{\partial}{\partial \theta}\widehat L_{n}(\theta) - \frac{\partial}{\partial \theta} L_{n}(\theta)  \Big\|_\Theta
         &\leq 
 \frac{1}{\sqrt{n}} \sum_{t=1}^{n} \Big\|\frac{\partial}{\partial \theta}\widehat{\ell}_{t}(\theta)-\frac{\partial}{\partial \theta}\ell_{t}(\theta) \Big\|_\Theta
 \leq 
 \frac{1}{\sqrt{n}}\sum_{i=1}^{m}\sum_{t=1}^{n}\Big\|\frac{\partial}{\partial \theta}\widehat{\ell}_{t,i}(\theta)-\frac{\partial}{\partial \theta}\ell_{t,i}(\theta) \Big\|_\Theta .  
 \end{align} 
Moreover, by proceeding as in Lemma 7.1 of \cite{Doukhan2015}, we can use  \textbf{A}$_i(\Theta)$ ($i=0,1$),  (\ref{moment}) and the condition (\ref{cond.th2}) to establish that
\[
 \E\Big[\frac{1}{\sqrt{n}} \sum_{t=1}^{n}\Big\|\frac{\partial}{\partial \theta}\widehat{\ell}_{t,i}(\theta)-\frac{\partial}{\partial \theta}\ell_{t,i}(\theta) \Big\|_\Theta\Big]   \limiten 0 \text{ for all } i=1,\cdots,m.
 \]
Thus, we can conclude the proof of (i) from (\ref{dif_Ln_hat.Ln}). \\
(ii) It holds that
\begin{align*}   
         \frac{1}{n} \Big\|\frac{\partial^2}{\partial \theta \partial \theta^T}\widehat L_n (\theta)- \frac{\partial^2}{\partial \theta \partial \theta^T} L_n (\theta)\Big\|_\Theta
         &\leq 
 \frac{1}{n} \sum_{t=1}^{n} \Big\|\frac{\partial^2}{\partial \theta \partial \theta^T}\widehat{\ell}_{t}(\theta)-\frac{\partial^2}{\partial \theta \partial \theta^T}\ell_{t}(\theta) \Big\|_\Theta \\
 & \leq 
 \frac{1}{n}\sum_{i=1}^{m}\sum_{t=1}^{n}\Big\|\frac{\partial^2}{\partial \theta \partial \theta^T}\widehat{\ell}_{t,i}(\theta)-\frac{\partial^2}{\partial \theta \partial \theta^T}\ell_{t,i}(\theta) \Big\|_\Theta .  
 \end{align*}
By going as in the proof of Lemma 7.1 of \cite{Doukhan2015}, one easily get for $i=1,\cdots,m$,
$ \frac{1}{n} \sum_{t=1}^{n}\Big\|\frac{\partial^2}{\partial \theta \partial \theta^T}\widehat{\ell}_{t,i}(\theta)-\frac{\partial^2}{\partial \theta \partial \theta^T}\ell_{t,i}(\theta) \Big\|_\Theta = o(1)$,
which shows that (ii) holds. 
\begin{flushright}
$\blacksquare$
\end{flushright}

\noindent 
 The following lemma will also be   needed.
 \begin{lem}\label{lem3} 
 If the assumptions of Theorem \ref{th2} hold, then
\begin{enumerate}
\item [(a)]  the matrices
$
J_{\theta^*} = \E\big[\frac{\partial \lambda^T_{0}(\theta^* )}{ \partial \theta} {D}^{-1}_0(\theta^*)   \frac{\partial \lambda_{0}(\theta^* )}{ \partial \theta^T} \big]  
       \it    ~~and~~ 
           I_{\theta^*}= \E \big[ \frac{\partial \lambda^T_{0}(\theta^* )}{ \partial \theta}  {D}^{-1}_0(\theta^*) \Gamma_0 (\theta^*)  {D}^{-1}_0(\theta^*) \frac{\partial \lambda_{0}(\theta^* )}{ \partial \theta^T} \big] 
          $ 
            exist and are positive definite;
            
  \item [(b)] 
   $\left(\frac{\partial \ell_t(\theta^*)}{\partial \theta},\mathcal{F}_{t}\right)_{t \in \mathbb{Z}}$ is a stationary ergodic, square integrable martingale difference sequence with covariance matrix $I_{\theta^*}$;
  
  \item [(c)] $\E\big[ \| \frac{\partial^{2} \ell_0 (\theta)}{\partial \theta\partial \theta^T} \|_\Theta \big] < \infty$ and    
    $\E\big[\frac{\partial^{2} \ell_0 (\theta^*)}{\partial \theta\partial \theta^T}\big]=-J_{\theta^*}$;

\item [(d)]  $J(\widehat \theta_{n}) \limitepsn J_{\theta^*}$ and that the matrix $J_{\theta^*}$ is invertible. 
\end{enumerate}
 \end{lem} 
 
 \medskip
 \noindent
 \emph{\bf Proof.}
 \begin{enumerate}
\item [($a$)]  From the assumption ($\mathcal{MOD}.\textbf{A0}$), we can find a constant $C>0$ such that, it holds a.s. 
  \[ \E \Big \| \frac{\partial \lambda^T_{0}(\theta )}{ \partial \theta}  {D}^{-1}_0(\theta)   \frac{\partial \lambda_{0}(\theta )}{ \partial \theta^T} \Big \|_\Theta 
  \leq  C \E \Big \| \frac{\partial \lambda_{0}(\theta )}{ \partial \theta} \Big \|^2_\Theta \leq   C \sum_{i=1}^m \E \Big \| \frac{\partial \lambda_{0,i}(\theta )}{ \partial \theta} \Big \|^2_\Theta .   \]
One can show as in the proof of Lemma 7.1 of \cite{Doukhan2015} that for any $i=1,\cdots,m$, 
$\E \big \| \frac{\partial \lambda_{0,i}(\theta )}{ \partial \theta} \big \|^2_\Theta < \infty$. Hence,
\[ \E \Big \| \frac{\partial \lambda^T_{0}(\theta )}{ \partial \theta}  {D}^{-1}_0(\theta)   \frac{\partial \lambda_{0}(\theta )}{ \partial \theta^T} \Big \|_\Theta < \infty ,  \] 
which establishes the existence of $J_{\theta^*}$.

\medskip
Using  ($\mathcal{MOD}.\textbf{A0}$) and H\"{o}lder's inequality, we have
\begin{align*}
&\E \Big \|  \frac{\partial \lambda^T_{0}(\theta )}{ \partial \theta}  {D}^{-1}_0(\theta) \Gamma_0 (\theta)  {D}^{-1}_0(\theta) \frac{\partial \lambda_{0}(\theta )}{ \partial \theta^T}
\Big \|_\Theta   
\leq 
   \E \Big[\left\| {D}^{-1}_0(\theta) (Y_0-\lambda_{0}(\theta ))\right\|^2_\Theta 
  \Big \| \frac{\partial \lambda_{0}(\theta )}{ \partial \theta} \Big \|^2_\Theta  \Big]\\ 
  &\leq 
 C \E \Big[ \big( \|Y_0 \|^2 + \|\lambda_0(\theta) \|^2_\Theta  \big) 
 \Big \| \frac{\partial \lambda_{0}(\theta )}{ \partial \theta} \Big \|^2_\Theta \Big] 
 \leq 
 C\left(  \E \| Y_0 \|^4 + \E\|\lambda_0(\theta) \|^4_\Theta \right  )^{1/2} 
 \Big(\E \Big \| \frac{\partial \lambda_{0}(\theta )}{ \partial \theta} \Big \|^4_\Theta\Big)^{1/2} \\
 &\leq 
 C \left(  \E \| Y_0 \|^4 + \sum_{i=1}^m \E\|\lambda_{0,i}(\theta) \|^4_\Theta \right  )^{1/2} 
 \Big(\sum_{i=1}^m \E \Big \| \frac{\partial \lambda_{0,i}(\theta )}{ \partial \theta} \Big \|^4_\Theta\Big)^{1/2}.
\end{align*}
According to the existence of the moment of order $4$ (from (\ref{moment}) with $\epsilon \geq 3$), $\E \| Y_0\|^4<\infty$. 
Furthermore, proceeding as in the proof of Theorem 3.1 and Lemma 7.1 in \cite{Doukhan2015}, one can also get 
$\E\|\lambda_{0,i}(\theta) \|^4_\Theta < \infty$ and
 $\E \big \| \frac{\partial \lambda_{0,i}(\theta )}{ \partial \theta} \big \|^4_\Theta <\infty$ for any $i=1,\cdots,m$. 
 Therefore,
\begin{equation}\label{I_theta_Theta}
\E \Big \|  \frac{\partial \lambda^T_{0}(\theta )}{ \partial \theta}  {D}^{-1}_0(\theta) \Gamma_0 (\theta)  {D}^{-1}_0(\theta) \frac{\partial \lambda_{0}(\theta )}{ \partial \theta^T}
\Big \|_\Theta < \infty,
\end{equation} 
  which establishes that $I_{\theta^*}$ exists. 

 \medskip
 \noindent
Now, let $U \in \R^d$ be a non zero vector. 
We have $ \frac{\partial \lambda_{0}(\theta^* )}{ \partial \theta^T} \cdot U \neq 0$ a.s. from the assumption ($\mathcal{MOD}.\textbf{A2}$);
 which implies
\[ U^T J_{\theta^*} U = \E\big[ U^T \frac{\partial \lambda^T_{0}(\theta^* )}{ \partial \theta}  {D}^{-1}_0(\theta^*)   \frac{\partial \lambda_{0}(\theta^* )}{ \partial \theta^T} U \big] 
= \E\big[ \big( {D}^{-1/2}_0(\theta^*) \frac{\partial \lambda_{0}(\theta^* )}{ \partial \theta^T} U \big)^T \cdot  \big( {D}^{-1/2}_0(\theta^*) \frac{\partial \lambda_{0}(\theta^* )}{ \partial \theta^T} U \big) \big]  > 0 \]
and
\begin{align*}
U^T I_{\theta^*} U
&=
\E\big[ U^T \frac{\partial \lambda^T_{0}(\theta^* )}{ \partial \theta}  {D}^{-1}_0(\theta^*) \Gamma_0 (\theta^*)  {D}^{-1}_0(\theta^*) \frac{\partial \lambda_{0}(\theta^* )}{ \partial \theta^T} U \big]\\
&=
 \E\big[ \big( \left(Y_0-\lambda_{0}(\theta^* )\right)^T {D}^{-1/2}_0(\theta^*) \frac{\partial \lambda_{0}(\theta^* )}{ \partial \theta^T} U \big)^T \cdot  \big( \left(Y_0-\lambda_{0}(\theta^* )\right)^T {D}^{-1/2}_0(\theta^*) \frac{\partial \lambda_{0}(\theta^* )}{ \partial \theta^T} U \big) \big] > 0 .
\end{align*}
Hence, $J_{\theta^*}$ and $I_{\theta^*}$ are positive definite.

\item [($b$)] For any $\theta \in \Theta$, we have
     \begin{equation}\label{Deriv1_ell_t}
     \frac{\partial \ell_t(\theta)}{\partial \theta} =
   \sum_{i=1}^{m}
    \big(\frac{Y_{t,i}}{\lambda_{t,i}(\theta)}-1\big)  \frac{\partial}{\partial \theta}\lambda_{t,i}(\theta)  
    = \frac{\partial \lambda^{T}_{t}(\theta)}{ \partial \theta} {  D}^{-1}_t(\theta) (Y_t- \lambda_{t}(\theta)). 
    \end{equation}
Then, according to the stability properties of $\{Y_{t},~t\in \Z\}$, the process $\big\{\frac{\partial \ell_t(\theta)}{\partial \theta},~t\in \Z \big\}$ is also stationary and ergodic.   
Moreover, since $\lambda_t(\theta^*)$ and $\frac{\partial \lambda_t(\theta^*)}{\partial \theta}  $ are $\mathcal{F}_{t-1}$-measurable, we have
\begin{align*}
\E \Big[ \frac{\partial \ell_t(\theta^*)}{\partial \theta}\Big]
&=
\sum_{i=1}^{m} \E\Big[ \E\big[\frac{\partial \lambda^{T}_{t}(\theta^*)}{ \partial \theta} {  D}^{-1}_t(\theta^*) (Y_t- \lambda_{t}(\theta^*))  | \mathcal{F}_{t-1}\big]\Big]\\
&=
\sum_{i=1}^{m} \E\Big[\frac{\partial \lambda^{T}_{t}(\theta^*)}{ \partial \theta} {  D}^{-1}_t(\theta^*)  \cdot \E\big[\left(Y_t- \lambda_{t}(\theta^*)\right)  | \mathcal{F}_{t-1}\big]\Big]=0
\end{align*}
In addition, 
$\E \big[ \frac{\partial \ell_0(\theta^*)}{\partial \theta}  \frac{\partial \ell_0(\theta^*)}{\partial \theta^{T}}\big] =I_{\theta^*}$. Hence, the part second part of the lemma holds. 

\item [($c$)] We have,
   \[ \E \Big \| \frac{\partial^2 \ell_0(\theta)}{\partial \theta \partial \theta^{T}} \Big \|_\Theta = \sum_{i=1}^m \E \Big \| \frac{\partial^2 \ell_{0,i}(\theta)}{\partial \theta \partial \theta^{T}} \Big \|_\Theta < \infty,\]
   where the above inequality holds since $\E \Big \| \frac{\partial^2 \ell_{0,i}(\theta)}{\partial \theta \partial \theta^{T}} \Big \|_\Theta < \infty$ for $i=1,\cdots,m$  by going as in the proof of Lemme 7.2 in \cite{Doukhan2015}. 
Moreover, according to (\ref{Deriv1_ell_t}), for any $\theta \in \Theta$, we have
     \begin{equation*}
     \frac{\partial^2 \ell_t(\theta)}{\partial \theta \partial \theta^{T}} =
     \sum_{i=1}^{m}
    \big(\frac{Y_{t,i}}{\lambda_{t,i}(\theta)}-1\big)  \frac{\partial^2 \lambda_{t,i}(\theta)}{\partial \theta \partial \theta^{T}}  
    -
     \sum_{i=1}^{m}
    \frac{Y_{t,i}}{\lambda^2_{t,i}(\theta)}  \frac{\partial \lambda_{t,i}(\theta)}{\partial \theta} \frac{\partial \lambda_{t,i}(\theta)}{\partial \theta^{T}}  .
    \end{equation*}
Then, using conditional expectations, we obtain
 \begin{equation*}
     \E\big[\frac{\partial^2 \ell_0(\theta^*)}{\partial \theta \partial \theta^{T}} \big]
     =
       -
     \E\big[ \sum_{i=1}^{m}
    \frac{1}{\lambda_{0,i}(\theta^*)}  \frac{\partial \lambda_{0,i}(\theta^*)}{\partial \theta} \frac{\partial \lambda_{0,i}(\theta^*)}{\partial \theta^{T}} \big]
    =
  - \big[\frac{\partial \lambda^{T}_{0}(\theta^* )}{ \partial \theta}  {D}^{-1}_0(\theta^*)   \frac{\partial \lambda_{0}(\theta^* )}{ \partial \theta^{T}} \big] 
   \it
   =- J_{\theta^*} . 
 \end{equation*}

\item [($d$)] We have 
\[
J(\widehat{\theta}_n) = \big(-\frac{1}{n}\frac{\partial^2}{\partial \theta \partial \theta_i} L_n (\bar \theta_{n,i})\big)_{1 \leq i \leq d}
=
\big(-\frac{1}{n}\sum_{t=1}^{n}\frac{\partial^{2}  \ell_t(\bar \theta_{n,i})}{\partial \theta\partial \theta_i}\big)_{1 \leq i \leq d}.
\]
Since $\widehat{\theta}_n  \limitepsn \theta^*$, $~ \bar \theta_{n,i} ~ \limitepsn \theta^*$ (for any  $i=1,\cdots,d$)
 and that $\E\big[\frac{\partial^{2} \ell_0 (\theta^*)}{\partial \theta\partial \theta^{T}}\big]=-J_{\theta^*}$ exists, by the uniform strong law of large numbers, for any $i =1,\cdots,d$, we get
\[ \frac{1}{n}\sum_{t=1}^{n}\frac{\partial^{2}  \ell_t(\bar \theta_{n,i})}{\partial \theta\partial \theta_i} 
\equalpsn
 \frac{1}{n}\sum_{t=1}^{n}\frac{\partial^{2}  \ell_t(\theta^*)}{\partial \theta\partial \theta_i}
\overset{a.s.}{\longrightarrow}  \E\big[\frac{\partial^{2} \ell_0(\theta^*)}{\partial \theta\partial\theta_i}\big] \text{ as } n \rightarrow \infty. 
 \]
Therefore,
\[
J(\widehat{\theta}_n) = \big(-\frac{1}{n}\sum_{t=1}^{n}\frac{\partial^{2}  \ell_t(\bar \theta_{n,i})}{\partial \theta\partial \theta_i}\big)_{1 \leq i \leq d}
\limitepsn  -\Big(\E\big[\frac{\partial^{2} \ell_0(\theta^*)}{\partial \theta\partial\theta_i}\big]\Big)_{1\leq i \leq d}=J_{\theta^*}.
\]
This completes the proof of the lemma.
\end{enumerate}
\begin{flushright}
$\blacksquare$ 
\end{flushright}

\medskip
 \noindent
Now, let us use the results of Lemmas \ref{lem2} and \ref{lem3} to complete the proof of Theorem \ref{th2}.\\
Since 
  $\widehat{\theta}_n$ is a local maximum of the function $\theta \mapsto  \widehat L_n (\theta)$  for $n$ large enough (from the assumption ($\mathcal{MOD}.\textbf{A1}$) and the consistency of $\widehat{\theta}_n$), $\frac{\partial}{\partial \theta}  \widehat L_n ( \widehat{\theta}_n)=0$.
  Thus, according to Lemma \ref{lem2}, the relation (\ref{mvt.1_dl}) becomes  
 \begin{equation}\label{mvt.2_dl}
\sqrt{n}J(\widehat \theta_{n}) (\widehat{\theta}_n-\theta^*) = 
\frac{1}{\sqrt{n}} \frac{\partial}{\partial \theta} L_n (\theta^*) + o_P(1).
\end{equation} 
Moreover, applying the central limit theorem to 
the sequence
 $\left(\frac{\partial \ell_t(\theta^*)}{\partial \theta} ,\mathcal{F}_{t}\right)_{t \in \mathbb{Z}}$, it holds that
\begin{equation*}
\frac{1}{\sqrt{n}} \frac{\partial }{\partial \theta} L_{n} (\theta^*) = \frac{1}{\sqrt{n}}\sum_{t=1}^{n}\frac{\partial}{\partial \theta} \ell_t(\theta^*) \limiteloin \mathcal{N}_d(0,I_{\theta^*}).
\end{equation*}
Therefore, for $n$ large enough, using Lemma \ref{lem3}(d) and the relation (\ref{mvt.2_dl}), we obtain
\[
\sqrt{n} (\widehat{\theta}_n-\theta^*) = 
J^{-1}_{\theta^*} \Big[ \frac{1}{\sqrt{n}} \frac{\partial}{\partial \theta} L_n (\theta^*)\Big] + o_P(1)
\limiteloin \mathcal{N}_d (0,J^{-1}_{\theta^*}I_{\theta^*} J^{-1}_{\theta^*}).
\]
This establishes the theorem.
\begin{flushright}
$\blacksquare$
\end{flushright}

\bibliographystyle{acm}

\end{document}